\documentclass[11pt]{amsart}
\usepackage[english]{babel}
\usepackage{amssymb,enumerate,bbm}

\newcommand\E{\mathbb{E}}
\newcommand\R{\mathbb{R}}
\newcommand\e{\varepsilon}

\numberwithin{equation}{section}

\def\diag{\operatorname*{diag}}

\catcode`\<=\active \def<{
\fontencoding{T1}\selectfont\symbol{60}\fontencoding{\encodingdefault}}
\catcode`\>=\active \def>{
\fontencoding{T1}\selectfont\symbol{62}\fontencoding{\encodingdefault}}
\catcode`\|=\active \def|{
\fontencoding{T1}\selectfont\symbol{124}\fontencoding{\encodingdefault}}
\newcommand{\assign}{:=}
\newcommand{\dueto}[1]{{\textbf{(#1) }}}
\newcommand{\mathd}{\mathrm{d}}
\newcommand{\nocomma}{}
\newcommand{\tmem}[1]{{\em #1\/}}

\newcommand{\tmmathbf}[1]{\ensuremath{\boldsymbol{#1}}}
\newcommand{\tmname}[1]{\textsc{#1}}
\newcommand{\tmop}[1]{\ensuremath{\operatorname{#1}}}
\newcommand{\tmstrong}[1]{\textbf{#1}}
\newcommand{\tmtextit}[1]{{\itshape{#1}}}
\newenvironment{enumeratenumeric}{\begin{enumerate}[1.] }{\end{enumerate}}

\newtheorem{theorem}{Theorem}[section]
\newtheorem*{theorem*}{Theorem}
\newtheorem{lemma}[theorem]{Lemma}
\newtheorem{cor}[theorem]{Corollary}
\theoremstyle{definition}

\newtheorem{dfn}[theorem]{Definition}
\renewcommand{\paragraph}[1]{ {\mbox{}\\ \noindent\textbf{#1.}} }

\title[ Semi-discrete Riesz transforms of the second order]{Various sharp estimates for semi-discrete Riesz transforms of the second order}

\author{K.\;Domelevo}

\author{A.\;Os\c ekowski}
	
\author{S.\;Petermichl}
\thanks{Adam Osekowski is supported by Narodowe Centrum Nauki Poland (NCN),
grant DEC-2014/14/E/ST1/00532}
\thanks{Stefanie Petermichl is supported by ERC project CHRiSHarMa DLV-862402}

\date{\today}
\begin{document}
{\maketitle}

\begin{abstract}
  We give several sharp estimates for a class of combinations of second order Riesz transforms on Lie
  groups $\mathbbm{G}=\mathbbm{G}_{x} \times \mathbbm{G}_{y}$ that are
  multiply connected, composed of a discrete abelian component
  $\mathbbm{G}_{x}$ and a connected component $\mathbbm{G}_{y}$ endowed with a biinvariant measure. 
  These estimates include new sharp $L^p$ estimates via Choi type constants, depending upon the multipliers of the operator. They also include weak-type, logarithmic and exponential estimates. We give an optimal $L^q \to L^p$ estimate as well.
 
It was shown recently by {\tmname{Arcozzi}}-{\tmname{Domelevo}}-{\tmname{Petermichl}} that such second order Riesz transforms applied to a function may be written as conditional expectation of a simple transformation of a stochastic integral associated with the function. 
  
 The proofs of our theorems combine this stochastic integral representation with a number of deep estimates for pairs of martingales under strong differential subordination by {\tmname{Choi}}, {\tmname{Banuelos}} and {\tmname{Os\c ekowski}}.
 
 When two continuous directions are available, sharpness is shown via the laminates technique.
 We show that sharpness is preserved in the discrete case  using Lax-Richtmyer theorem.

\end{abstract}


\section{Introduction}

Sharp, classical $L^{p}$ norm inequalities for pairs of differentially subordinate martingales
date back to the celebrated work of {\tmname{Burkholder}} {\cite{Bur1984a}} in
1984 where the optimal constant is exhibited. See also from the same author
{\cite{Bur1988a}}{\cite{Bur1991a}}. The relation between differentially
subordinate martingales and CZ (i.e. Cald{\'e}ron--Zygmund) operators is known at least since
{\tmname{Gundy}}--{\tmname{Varopoulos}} {\cite{GunVar1979}}.
{\tmname{Banuelos}}--{\tmname{Wang}} {\cite{BanWan1995}} were the first to
exploit this connection to prove new sharp inequalities for singular
intergrals. This intersection of probability theory with classical questions in harmonic analysis has lead 
to much interest and a vast literature has been accumulating on this line of
research.

\medskip

In this article we state a number of sharp estimates that hold in the very recent, new direction concerning the  
semi--discrete setting, applying it to a family of second order Riesz
transforms on multiply--connected Lie groups.  We recall their representation 
through stochastic integrals using jump processes on multiply--connected Lie groups from \cite{ArcDomPet2016a}. In this representation formula jump processes play a role, but the strong differential subordination holds between the martingales representing the test function and the operator applied to the test function.

\medskip

The usual procedure for obtaining (sharp) inequalities for operators of Calder\'on--Zygmund type from
inequalities for martingales is the following. Starting with a test
function $f$, martingales are built using Brownian motion or background noise and harmonic functions in the upper half space $\mathbbm{R}^{+} \times \mathbbm{R}^{n}$. Through the use of It\o formula, it is shown that the martingale arising in this way from $R f$, where $R$ is a Riesz
transform in $\mathbbm{R}^{n}$, is a martingale transform of the martingale arising from 
$f$. The two form a pair of martingales that have differential subordination and (in case of Hilbert or Riesz transforms)
orthogonality.  One then derives sharp
martingale inequalities under hypotheses of strong differential subordination
(and orthogonality) relations.

In the case of Riesz transforms of the second order, the use of heat extensions in
the upper half space instead of Poisson extensions originated in the context of a weighted estimate in
{\tmname{Petermichl}}--{\tmname{Volberg}} {\cite{PetVol2002a}} and was used to
prove $L^{p}$ estimates for the second order Riesz transforms based on the
results of Burkholder in {\tmname{Nazarov}}--{\tmname{Volberg}}
{\cite{VolNaz2004a}} as part of their best-at-time estimate for the
Beurling--Ahlfors operator, whose real and imaginary parts themselves are
second order Riesz transforms. We mention the recent version on discrete abelian groups
{\tmname{Domelevo--Petermichl}} {\cite{DomPet2014c}} also using a type of heat flow. These proofs are deterministic. The
technique of Bellman functions was used.
This deterministic strategy does well when no orthogonality is present and when strong subordination is the only important property. Stochastic proofs (aside from giving better estimates in some situations) also have the advantage that once the integral representation is known, the proofs are a very concise consequence of the respective statements on martingales. 

\medskip

In \cite{ArcDomPet2016a} the authors proved sharp $L^p$ estimates for semi--discrete
second order Riesz transforms $R_{\alpha}^2$ using stochastic integrals. There is an array of Riesz transforms of the second order that are treated, indexed my a matrix index $\alpha$ (see below for precisions on acceptable $\alpha$).
The following representation formula of semi-discrete second order Riesz transforms
$R^{2}_{\alpha}$ \tmtextit{{\`a} la} {\tmname{Gundy--Varopoulos}} (see
{\cite{GunVar1979}}) is instrumental:

\begin{theorem*}
\dueto{Arcozzi--Domelevo--Petermichl, 2016}
\label{T: a la Gundy-Varopoulos}
The second order Riesz transform
  $R_{\alpha}^{2} f$ of a function $f \in L^{2} ( \mathbbm{G} )$ as defined in
  (\ref{eq: definition 2nd order Riesz transforms}) can be written as the
  conditional expectation $$\mathbbm{E} ( M_{0}^{\alpha ,f} | \mathcal{Z}_{0}
  =z ).$$ Here $M_{t}^{\alpha ,f}$ is a suitable martingale transform of a
  martingale $M^{f}_{t}$ associated to $f$, and $\mathcal{Z}_{t}$ is a
  suitable random walk on $\mathbbm{G}$
\end{theorem*}

We remark that the $L^{p}$ estimates of the discrete Hilbert transform on the integers are still
open. It is a famous conjecture that this operator has the same norm as its
continuous counterpart. 

\medskip

These known $L^{p}$ norm inequalities use special functions found in
the results of
{\tmname{Pichorides}} {\cite{Pic1972}}, \tmname{Verbitsky}\cite{Ver1980a}, {\tmname{Ess{\'e}n}} {\cite{Ess1984}},
{\tmname{Banuelos}}--{\tmname{Wang}} {\cite{BanWan1995}} when orthogonality is
present in addition to differential subordination or {\tmname{Burkholder}}
{\cite{Bur1984a}}{\cite{Bur1987a}}{\cite{Bur1988a}}, {\tmname{Wang}}
{\cite{Wan1995a}} when differential subordination is the only hypothesis.

\medskip

The aim of the present paper is to establish new estimates for semi--discrete Riesz
transforms by using the martingale representation above together
with recent martingale inequalities found in the literature.

\medskip

Here is a brief description of the new results in this paper. 

\begin{itemize}
\item
In the case where the function $f$ is real valued, we can obtain better estimates for $R_{\alpha}^2$ than in the general case. These estimates depend upon the make of the matrix index $\alpha$. The precise statement is found in Theorem \ref{T: Choi constant estimate}. 

\item
We prove a refined sharp weak type estimate using a weak type norm defined just before the statement of Theorem
  \ref{T: weak-type estimate} .

\item
We prove logarithmic and exponential estimates, in a sense limiting (in $p$) cases of the classical sharp $L^p$ estimate. See Theorem \ref{T: logarithmic and exponential}.

\item
We consider the norm estimates of the $R^2_{\alpha}:L^q\to L^p$, spaces of different exponent. The statement is found in Theorem  \ref{T: mixed estimate}.

\end{itemize}

\subsection{Differential operators and Riesz transforms}

\mbox{}\smallskip

\paragraph{First order derivatives and tangent planes}
We will consider Lie groups $\mathbbm{G} \assign \mathbbm{G}_{x} \times
\mathbbm{G}_{y} ,$ where $\mathbbm{G}_{x}$ is a discrete abelian group with a
fixed set $G$ of $m$ generators, and their reciprocals, and $\mathbbm{G}_{y}$
is a connected, Lie group of dimension $n$ endowed with a biinvariant
metric. The choice of the set $G$ of generators in $\mathbbm{G}_{x}$
corresponds to the choice of a bi-invariant metric structure on
$\mathbbm{G}_{x}$. We will use on $\mathbbm{G}_{x}$ the multiplicative
notation for the group operation. We will define a product metric structure on
$\mathbbm{G}$, which agrees with the Riemannian structure on the first factor,
and with the discrete ``word distance'' on the second. We will at the same
time define a ``tangent space'' \ $T_{z} \mathbbm{G}$ for $\mathbbm{G}$ at a
point $z= ( x,y ) \in ( \mathbbm{G}_{x} \times \mathbbm{G}_{y} )
=\mathbbm{G}$. We will do this in three steps.

First, since $\mathbbm{G}_{y}$ is an $n$-dimensional connected Lie group with
Lie algebra $\mathfrak{G}_{y}$. We can identify each left-invariant vector
field $Y$ in $\mathfrak{G}_{y}$ with its value at the identity $e$,
$\mathfrak{G}_{y} \equiv T_{e} \mathbbm{G}_{y}$. Since $\mathbbm{G}$ is
compact, it admits a bi-invariant Riemannian metric, which is unique up to a
multiplicative factor. We normalize it so that the measure $\mu_{y}$
associated with the metric satisfies $\mu_{y} ( \mathbbm{G}_{y} ) =1$. The
measure $\mu_{y}$ is also the normalized Haar measure of the group. We denote
by $< \cdot , \cdot >_{y}$ be the corresponding inner product on $T_{y}
\mathbbm{G}_{y}$ and by $\tmmathbf{\nabla_{y}} f ( y )$ the gradient at $y \in
\mathbbm{G}_{y}$ of a smooth function $f:\mathbbm{G}_{y} \rightarrow
\mathbbm{R}$. Let $Y_{1} , \ldots ,Y_{n}$ be a orthonormal basis for
$\mathfrak{G}_{y}$. The gradient of $f$ can be written $\tmmathbf{\nabla_{y}}
f=Y_{1} ( f ) Y_{1} + \ldots +Y_{n} ( f ) Y_{n}$.

Second, in the discrete component $\mathbbm{G}_{x}$, let $\mathfrak{G}_{x} = (
g_{i} )_{i=1, \ldots ,m}$ be a set of generators for $\mathbbm{G}_{x}$, such
that for $i \neq j$ and $\sigma = \pm 1$ we have $g_{i} \neq g_{j}^{\sigma}$.
The choice of a particular set of generators induces a word metric, hence, a
geometry, on $\mathbbm{G}_{x}$. Any two sets of generators induce bi-Lipschitz
equivalent metrics.

At any point $x \in \mathbbm{G}_{x}$, and given a direction $i \in \{ 1,
\ldots ,m \}$, we can define the right and the left derivative at $x$ in the
direction $i$:
\[ ( \partial^{+} f/ \partial x_{i} ) ( x,y ) \assign f ( x+ g_{i} ,y ) -f ( x
   ,y ) \assign ( \partial_{i}^{+} f ) ( x,y ) \]
\[ \hspace{1em} ( \partial^{-} f/ \partial x_{i} ) ( x,y ) \assign f ( x ,y )
   -f ( x-g_{i}  ,y ) \assign ( \partial_{i}^{-} f ) ( x,y ) . \]
Comparing with the continuous component, this suggests that the tangent plane
$\hat{T}_{x} \mathbbm{G}_{x}$ at a point $x$ of the discrete group
$\mathbbm{G}_{x}$ might actually be split into a ``right'' tangent plane
$T_{x}^{+} \mathbbm{G}_{x}$ and a ``left'' tangent plane $T_{x}^{-}
\mathbbm{G}_{x}$, according to the direction with respect to which discrete
differences are computed. We consequently define the {\bf{augmented}}
discrete gradient $\widehat{\tmmathbf{\nabla}}_{x} f ( x )$, with a
{\tmem{hat}}, as the $2m$--vector of $\hat{T}_{x} \mathbbm{G}_{x} \assign
T_{x}^{+} \mathbbm{G}_{x} \oplus T_{x}^{-} \mathbbm{G}_{x}$ accounting for all
the local variations of the function $f$ in the direct vicinity of $x$; that
is, the $2m$--column--vector
\[ \widehat{\tmmathbf{\nabla}}_{x} f ( x ) \assign ( X_{1}^{+} f,X_{2}^{+} f,
   \ldots ,X_{1}^{-} f,X_{2}^{-} f, \ldots ) ( x ) = \sum_{i=1}^{m} \sum_{\tau
   = \pm} X_{i}^{\tau} f ( x ) \hspace{1em} \] 
   with $X_{i}^{\tau}  \in \hat{T}_{x} \mathbbm{G}_{x}, $
where we noted the discrete derivatives $X_{i}^{\pm} f \assign
\partial_{i}^{\pm} f$ and introduced the discrete $2m$--vectors
$X_{i}^{\pm}$ as the column vectors of $\mathbbm{Z}^{2m}$
\[ X_{i}^{+} = ( 0, \ldots ,1, \ldots ,0 ) \times \tmmathbf{0}_{m} ,
   \hspace{1em} X_{i}^{-} =\tmmathbf{0}_{m} \times ( 0, \ldots ,1, \ldots ,0 )
   .\]
Here the $1$'s in $X_{i}^{\pm}$ are located at the $i$--th position of
respectively the first or the second $m$--tuple. Notice that those vectors are
independent of the point $x$. The scalar product on $\hat{T}_{x}
\mathbbm{G}_{x} \assign T_{x}^{+} \mathbbm{G}_{x} \oplus T_{x}^{-}
\mathbbm{G}_{x}$ is defined as
\[ ( U,V )_{\hat{T}_{x} \mathbbm{G}_{x}} \assign \frac{1}{2} \sum_{i=1}^{m}
   \sum_{\tau = \pm} U_{i}^{\tau} V^{\tau}_{i}  . \]
We chose to put a factor $\frac{1}{2}$ in front of the scalar product to
compensate for the fact that we consider both left and right differences.

Finally, for a function $f$ defined on the cartesian product $\mathbbm{G}
\assign \mathbbm{G}_{x} \times \mathbbm{G}_{y}$, the (augmented) gradient
$\widehat{\tmmathbf{\nabla}}_{z} f ( z )$ at the point $z= ( x,y )$ is an
element of the tangent plane $\hat{T}_{z} \mathbbm{G} \assign \hat{T}_{x}
\mathbbm{G}_{x} \oplus T_{y} \mathbbm{G}_{y}$, that is a $( 2m+n
)$--column--vector
\begin{eqnarray*}
  \widehat{\tmmathbf{\nabla}}_{z} f ( z ) & \assign & \sum_{i=1}^{m}
  \sum_{\tau = \pm} X_{i}^{\tau} f ( z )  \hat{X}_{i}^{\tau} +
  \sum_{j=1}^{n} Y_{j} f ( z )  \hat{Y}_{j} ( z )\\
  & = & ( X_{1}^{+} f,X_{2}^{+} f, \ldots ,X_{1}^{-} f,X_{2}^{-} f, \ldots
  ,Y_{1} f,Y_{2} f, \ldots ) ( z )
\end{eqnarray*}
where $\hat{X}_{i}^{\tau}$ and $\hat{Y}_{j} ( z )$ can be identified with
column vectors of size $( 2m+n )$ with obvious definitions and scalar product
$( \cdot , \cdot )_{\hat{T}_{z} \mathbbm{G}_{z}}$.

Let $\mathd \mu_{z} \assign \mathd \mu_{x} \mathd \mu_{y}$, $\mathd \mu_{x}$
being the counting measure on $\mathbbm{G}_{x}$ and $\mathd \mu_{y}$ being \
the Haar measure on $\mathbbm{G}_{y}$. The inner product of $\varphi , \psi$
in $L^{2} ( \mathbbm{G} )$ is
\begin{eqnarray*}
  ( \varphi , \psi )_{L^{2} ( \mathbbm{G} )} & \assign & \int_{\mathbbm{G}}
  \varphi ( z ) \psi ( z ) d \mu_{z} ( z ) .
\end{eqnarray*}

Finally, we make the following hypotheses
 
\paragraph{Hypothesis}We assume everywhere in the sequel:
\begin{enumeratenumeric}
  \item The discrete component $\mathbbm{G}_{x}$ of the Lie group
  $\mathbbm{G}$ is an abelian group
  
  \item The connected component $\mathbbm{G}_{y}$ of the Lie group
  $\mathbbm{G}$ is a Lie group that can be endowed with a biinvariant Riemannian
  metric, so that the family $( Y_{j} )_{j=1, \ldots ,n}$ commutes with
  $\Delta_{y}$.
\end{enumeratenumeric}
\mbox{}\\

Notice that this includes compact Lie groups  $\mathbbm{G}_{y}$ since those can be endowed with a biinvariant metric.
It also includes the usual Euclidian spaces since those are commutative.

\paragraph{Riesz transforms} Following {\cite{Arc1995a}}{\cite{Arc1998a}},
recall first that for a compact Riemannian manifold $\mathbbm{M}$ without
boundary, one denotes by $\tmmathbf{\nabla}_{\mathbbm{M}}$,
$\tmop{div}_{\mathbbm{M}}$ and $\Delta_{\mathbbm{M}} \assign
\tmop{div}_{\mathbbm{M}} \tmmathbf{\nabla}_{\mathbbm{M}}$ respectively the
gradient, the divergence and the Laplacian associated with $\mathbbm{M}$. Then
$- \Delta_{\mathbbm{M}}$ is a positive operator and the vector Riesz transform
is defined as the linear operator
\[ \tmmathbf{R}_{\mathbbm{M}} \assign \tmmathbf{\nabla}_{\mathbbm{M}} \circ (
   - \Delta_{\mathbbm{M}} )^{-1/2} \]
acting on $L^{2}_{0} ( \mathbbm{M} )$ ($L^{2}$ functions with vanishing mean).
It follows that if $f$ is a function defined on $\mathbbm{M}$ and $y \in
\mathbbm{M}$ then $\tmmathbf{R}_{\mathbbm{M}} f ( y )$ is a vector of the
tangent plane $T_{y} \mathbbm{M}$.

Similarly on $\mathbbm{M}=\mathbbm{G}$, we define
$\tmmathbf{\nabla}_{\mathbbm{G}} \assign \widehat{\tmmathbf{\nabla}}_{z}$ as
before, and then we define the divergence operator as its formal adjoint, that
is $- \tmop{div}_{\mathbbm{G}} =- \widehat{  \tmop{div}}_{z} \assign
\widehat{\tmmathbf{\nabla}}_{z}^{\ast}$, with respect to the natural $L^{2}$
inner product of vector fields:
\begin{eqnarray*}
  ( U,V )_{L^{2} ( \hat{T} \mathbbm{G} )} & \assign & \int_{\mathbbm{G}} ( U (
  z ) ,V ( z ) )_{\hat{T}_{z} \mathbbm{G}} \; \mathd \mu_{z} ( z )
\end{eqnarray*}
We have the $L^{2}$-adjoints $( X_{i}^{\pm} )^{\ast} =-X_{i}^{\mp}$ and
$Y_{j}^{\ast} =-Y_{j}$. If $U \in \hat{T} \mathbbm{G}$ is defined by
\[ U ( z ) = \sum_{i=1}^{m} \sum_{\tau = \pm} U_{i^{}}^{\tau} ( z )
   \hat{X}_{i}^{\tau} + \sum_{j=1}^{n} U_{j} ( z ) \hat{Y}_{j} , \]
we define its divergence $\widehat{\tmmathbf{\nabla}}_{z}^{\ast} U$ as
\begin{eqnarray*}
  \widehat{\tmmathbf{\nabla}}_{z}^{\ast} U ( z ) & \assign & - \frac{1}{2}
  \sum_{i=1}^{m} \sum_{\tau = \pm} X_{i}^{- \tau} U_{i}^{\tau} ( z ) -
  \sum_{j=1}^{n} Y_{j} U_{j} ( z ) .
\end{eqnarray*}
The Laplacian $\Delta_{\mathbbm{G}}$ is as one might expect:
\begin{eqnarray*}
  \Delta_{z} f ( z ) & \assign & - \widehat{\tmmathbf{\nabla}}_{z}^{\ast}
  \widehat{\tmmathbf{\nabla}}_{z} f ( z ) =-
  \widehat{\tmmathbf{\nabla}}_{x}^{\ast} \widehat{\tmmathbf{\nabla}}_{x} f ( z
  ) - \widehat{\tmmathbf{\nabla}}^{\ast}_{y} \widehat{\tmmathbf{\nabla}}_{y} f
  ( z )\\
  & = & \sum_{i=1}^{m} X_{i}^{-} X_{i}^{+} f ( z )  +
   \sum_{j=1}^{n} Y^{2}_{j} f ( z )\\
  & = & \sum_{i=1}^{m} X_{i}^{2} f ( z )  + 
  \sum_{j=1}^{n} Y^{2}_{j} f ( z )\\
  & =: & \Delta_{x} f ( z )  +  \Delta_{y} f ( z )
\end{eqnarray*}
where we denoted $X_{i}^{2} \assign X_{i}^{+} X_{i}^{-} =X_{i}^{-} X_{i}^{+}$.
We have chosen signs so that $- \Delta_{\mathbbm{G}} \geqslant 0$ as an
operator. The Riesz vector $( \hat{\tmmathbf{R}}_{z} f ) ( z )$ is the $( 2m+n
)$--column--vector of the tangent plane $\hat{T}_{z} \mathbbm{G}$ defined as
the linear operator
\[ \hat{\tmmathbf{R}}_{z} f \assign \left( \widehat{\tmmathbf{\nabla}}_{z} f
   \right) \circ ( - \Delta_{z} f )^{-1/2} \]
We also define transforms along the coordinate directions:
\[ R^{\pm}_{i} = X^{\pm}_{i} \circ ( - \Delta_{z} )^{-1/2} \hspace{1em}
   \tmop{and} \hspace{1em} R_{j} = Y_{j} \circ ( - \Delta_{z} )^{-1/2}
   . \]

\paragraph{Plan of the paper} In the next two sections, we present successively the main results of the paper
and recall the weak formulations involving second order Riesz transforms and semi-discrete heat extensions.
Section \ref{S: stochastic} introduces the stochastic setting for our problems. This includes in Subsection
\ref{SS: stochastic integrals and martingale transforms} semi-discrete
random walks, martingale transforms and quadratic covariations. Subsection \ref{SS: martingale inequalities} presents a set of martingale
inequalities already known in the literature. Finally, in Section \ref{S: proofs of the main results} we give the proof of the main results.

\subsection{Main results}
\label{SS: main results}

In this text, we are concerned with second order Riesz transforms and
combinations thereof. We first define the square Riesz transform in the
(discrete) direction $i$ to be
\[ R_{i}^{2} \assign R_{i}^{+} R_{i}^{-} =R_{i}^{-} R_{i}^{+} . \]
Then, given $\alpha \assign ( ( \alpha^{x}_{i} )_{i=1 \ldots m} , (
\alpha^{y}_{j k} )_{j,k=1 \ldots n} ) \in \mathbbm{C}^{m} \times
\mathbbm{C}^{n \times n}$, we define $R_{\alpha}^{2}$ to be the following
combination of second order Riesz transforms:
\begin{equation}
  R_{\alpha}^{2} \assign \sum^{m}_{i=1} \alpha^{x}_{i} \; R^{2}_{i} +
  \sum^{n}_{j,k=1} \alpha^{y}_{j k} \; R_{j} R_{k} , \label{eq: definition 2nd
  order Riesz transforms}
\end{equation}
where the first sum involves squares of discrete Riesz transforms as defined
above, and the second sum involves products of continuous Riesz transforms.
This combination is written in a condensed manner as the quadratic form
\[ R^{2}_{\alpha} = \left( \hat{\tmmathbf{R}}_{z} , \tmmathbf{A}_{\alpha}
   \hat{\tmmathbf{R}}_{z} \right) \]
where $\tmmathbf{A}_{\alpha}$ is the $( 2m+n ) \times ( 2m+n )$ block matrix
\begin{equation}
  \tmmathbf{A}_{\alpha} \assign \left(\begin{array}{cc}
    \tmmathbf{A}_{\alpha}^{x} & \tmmathbf{0}\\
    \tmmathbf{0} & \tmmathbf{A}_{\alpha}^{y}
  \end{array}\right) \label{eq: Matrix Aalpha}
\end{equation}
with
\[ \tmmathbf{A}_{\alpha}^{x} = \tmmathbf{\tmop{diag}} ( \alpha_{1}^{x} ,
   \ldots , \alpha_{m}^{x} , \alpha_{1}^{x} , \ldots , \alpha_{m}^{x} )
   \in \mathbbm{C}^{2m \times 2m}, 
   \tmmathbf{A}_{\alpha}^{y} = ( \alpha^{y}_{j k} )_{j,k=1 \ldots n}
    \in \mathbbm{C}^{n \times n} . \]

In the theorems below, we assume that $\mathbbm{G}$ is a Lie group and
  $R^{2}_{\alpha}$ 
 is a combination of second order Riesz transforms as defined above. The first application of the stochastic integral formula, Theorem \ref{T: LpRieszsemidiscrete} was done in \cite{ArcDomPet2016a}, while the other applications, Theorems \ref{T: Choi constant estimate}  \ref{T: weak-type estimate}  \ref{T: logarithmic and exponential} and  \ref{T: mixed estimate} are new.

\begin{theorem}\label{T: LpRieszsemidiscrete}{\dueto{Arcozzi--Domelevo--Petermichl,
  2016}}
  \label{T: p minus 1 estimate} For any $1<p<\infty$ we have
  \[ \| R_{\alpha}^{2} \|_p \leqslant \left\| \tmmathbf{A}_{\alpha} \right\|_{2}
     \  ( p^{\ast} -1 ) , \]
  where, as previously, $p^*=\max\{p,p/(p-1)\}$. 
\end{theorem}

Above, we have set: $$\left\| \tmmathbf{A}_{\alpha} \right\|_{2} = \max \left(
\left\| \tmmathbf{A}_{\alpha}^{x} \right\|_{2} , \left\|
\tmmathbf{A}_{\alpha}^{y} \right\|_{2} \right) = \max \left( | \alpha_{1}^{x}
| , \ldots , | \alpha_{m}^{x} | , \left\| \tmmathbf{A}_{\alpha}^{y}
\right\|_{2} \right).$$
In the case where $\mathbbm{G}=\mathbbm{G}_{x}$ only consists of the discrete
component, this was proved in {\cite{DomPet2014b}}{\cite{DomPet2014c}} using
the deterministic Bellman function technique. In the case where
$\mathbbm{G}=\mathbbm{G}_{y}$ is a connected compact Lie group, this was
proved in {\cite{BanBau2013}} using Brownian motions defined on manifolds and
projections of martingale transforms.

In the case where the function $f$ is real valued, we can obtain better estimates. For any real numbers $a<b$ and any $1<p<\infty$, let $\mathfrak{C}_{a,b,p}$ be the constants introduced in {\tmname{Ba\~nuelos}} and {\tmname{Os\c ekowski}} {\cite{BanOse2012a}}. 

\begin{theorem}
  \label{T: Choi constant estimate}Assume that $a\tmmathbf{I} \leqslant
  \tmmathbf{A}_{\alpha} \leqslant b\tmmathbf{I}$ in the sense of quadratic
  forms. Then $R^{2}_{\alpha} :L^{p} (
  \mathbbm{G},\mathbbm{R} ) \rightarrow L^{p} ( \mathbbm{G},\mathbbm{R} )$
  enjoys the norm estimate $\| R_{\alpha}^{2} \|_{p} \leqslant 
  \mathfrak{C}_{a, \nocomma b,p} \nocomma$.
\end{theorem}

We should point out here that the constants $\mathfrak{C}_{a,b,p}$ appear in earlier works of {\tmname{Burkholder}} \cite{Bur1984a} (for $a=-b$: then $\mathfrak{C}_{a,b,p}=b(p^\ast-1)$), and in the paper {\cite{Cho1992a}} by {\tmname{Choi}} (in the case when one of $a$, $b$ is zero). 
The Choi constants are not explicit; an approximation of
$\mathfrak{C}_{0,1,p}$ is known and writes as
\[ \mathfrak{C}_{0,1,p} = \tfrac{p}{2} + \tfrac{1}{2} \log \left(
   \tfrac{1+e^{-2}}{2} \right) + \tfrac{\beta_{2}}{p} + \ldots ., \]
    with $\beta_{2} = \log^{2} \left( \tfrac{1+e^{-2}}{2}
   \right) + \tfrac{1}{2} \log \left( \tfrac{1+e^{-2}}{2} \right) -2 \left(
   \tfrac{e^{-2}}{1+e^{-2}} \right)^{2} . $

Coming back to complex-valued functions, we will also establish the following weak-type bounds. We consider the norms
$$ ||f||_{L^{p,\infty}(\mathbbm{G},\mathbbm{C})}=\sup \left\{\mu_z(E)^{1/p-1}\int_{E}fd\mu_z\right\},$$
where the supremum is taken over the class of all measurable subsets $E$ of $\mathbbm{G}$ of positive measure.

\begin{theorem}
  \label{T: weak-type estimate} For any $1<p<\infty$ we have
  \[ \| R_{\alpha}^{2} \|_{L^{p} ( \mathbbm{G},\mathbbm{C} ) \rightarrow L^{p,\infty} (
  \mathbbm{G},\mathbbm{C} )} \leqslant \left\| \tmmathbf{A}_{\alpha} \right\|_{2}\cdot
     \begin{cases}
\displaystyle \left(\frac{1}{2}\Gamma\left(\frac{2p-1}{p-1}\right)\right)^{1-1/p} & \mbox{if }1<p\leq 2,\\
\displaystyle \left(\frac{p^{p-1}}{2}\right)^{1/p} & \mbox{if }p\geq 2.
\end{cases}  \]
\end{theorem}

We will also prove the following logarithmic and exponential estimates, which can be regarded as versions of Theorem \ref{T: p minus 1 estimate} for $p=1$ and $p=\infty$. Consider the Young functions $\Phi,\,\Psi:[0,\infty)\to [0,\infty)$, given by 
$\Phi(t)=e^t-1-t$ and $\Psi(t)=(t+1)\log(t+1)-t.$

\begin{theorem}
  \label{T: logarithmic and exponential} Let $K>1$ be fixed.

(i) For any measurable subset $E$ of $\mathbbm{G}$ and any $f$ on $\mathbbm{G}$ we have
  \[ \int_E | R_{\alpha}^{2}f | d\mu_z \leqslant \left\| \tmmathbf{A}_{\alpha} \right\|_{2}\cdot \left(K\int_{\mathbbm{G}} \Psi(|f|)d\mu_z+\frac{\mu_z(E)}{2(K-1)}\right). \]
  
(ii) For any $f:\mathbbm{G}\to \mathbbm{C}$ bounded by $1$,
  \[ \int_{\mathbbm{G}} \Phi\left(\frac{| R_{\alpha}^{2}f |}{K\left\| \tmmathbf{A}_{\alpha} \right\|_{2}}\right) d\mu_z \leqslant \frac{||f||_{L^1(\mathbbm{G},\mathbbm{C})}}{2K(K-1)}. \]
\end{theorem}

Our final result concerns another extension of Theorem \ref{T: p minus 1 estimate}, which studies the action of $R_\alpha^2$ between two different $L^p$ spaces. For $1\leq p<q<\infty$, let $C_{p,q}$ be the constant defined by {\tmname{Os\c ekowski}} in \cite{Ose2010a}.

\begin{theorem}
  \label{T: mixed estimate} For any $1\leq p<q<\infty$, any measurable subset $E$ of $\mathbbm{G}$ and any $f\in L^q(\mathbbm{G})$ we have
  \[ \| R_{\alpha}^{2}f\|_{L^p(E,\mathbbm{C})} \leqslant C_{p,q}\left\| \tmmathbf{A}_{\alpha} \right\|_{2}\|f\|_{L^q(\mathbbm{G},\mathbbm{C})}\mu_z(E)^{1/p-1/q}
     . \]
\end{theorem}

An interesting feature is that all the estimates in the five theorems above are sharp when the group $\mathbbm{G}=\mathbbm{G}_{x}
  \times \mathbbm{G}_{y}$ and $\dim ( \mathbbm{G}_{y} ) + \dim^{\infty} (
  \mathbbm{G}_{x} ) \geqslant 2$, where $\dim^{\infty} ( \mathbbm{G}_{x} )$
  denotes the number of infinite components of $\mathbbm{G}_{x}$.

\subsection{Weak formulations}\label{SS: Weak formulations}

Let $f:\mathbbm{G} \rightarrow \mathbbm{C}$ be given. The heat extension
$\tilde{f} ( t )$ of $f$ is defined as $\tilde{f} ( t ) \assign e^{t
\Delta_{z}} f=:P_{t} f$. We have therefore $\tilde{f} ( 0 ) =f$. The aim of
this section is to derive weak formulations for second order Riesz transforms.
We start with the weak formulation of the identity operator $\mathcal{I}$,
that is obtained by using semi-discrete heat extensions (see \cite{ArcDomPet2016a} for details).

Assume $f$ in $L^{2} ( \mathbbm{G} )$ and $g$ in $L^{2} (  \mathbbm{G} )$. Let $\bar{f}$ be the average of $f$ on $\mathbbm{G}$
if $\mathbbm{G}$ has finite measure and zero otherwise. Then
  \begin{eqnarray*}
    ( \mathcal{I} f,g ) & = & ( f,g )_{L^{2} ( \mathbbm{G} )}\\
    & = & \bar{f}\bar{g} + 2 \int^{\infty}_{0} \left( \widehat{\tmmathbf{\nabla}}_{z} P_{t} f,
    \widehat{\tmmathbf{\nabla}}_{z} P_{t} g \right)_{L^{2} ( \hat{T}
    \mathbbm{G} )} \; \mathd  t\\
    & = &  \bar{f}\bar{g} + 2 \int^{\infty}_{0} \int_{z \in \mathbbm{G}} \left\{ \frac{1}{2}
    \sum_{i=1}^{m} \sum_{\tau = \pm} ( X_{i}^{\tau} P_{t} f ) ( z )
     ( X_{i}^{\tau} P_{t} g ) ( z ) \  + \right.\\
    &  & \left. \sum_{j=1}^{n} ( Y_{j} P_{t} f ) ( z ) ( Y_{j}
    P_{t} g ) ( z ) \right\}  \mathd \mu_{z} ( z )   \mathd  t
  \end{eqnarray*}
  and the sums and integrals that arise converge absolutely.

In order to pass to the weak formulation for the squares of Riesz transforms,
we first observe that the following commutation relations hold
\begin{eqnarray*}
    Y_{j} \circ \Delta_{z} & = & \Delta_{z} \circ Y_{j}\\
    X_{i}^{\tau} \circ \Delta_{z} & = & \Delta_{z} \circ X_{i}^{\tau} ,
    \hspace{1em} \tau \in \{ +,- \}
\end{eqnarray*}

This is an easy consequence of the hypothesis made on the Lie group.
Following \cite{ArcDomPet2016a}, the following weak formulation for second order Riesz transforms holds

Assume $f$ in $L^{2} ( \mathbbm{G} )$ and $g$ in $L^{2} (  \mathbbm{G} )$, then
  \begin{eqnarray*}
    ( R^{2}_{\alpha} f,g )_{L^{2} ( \mathbbm{G} )} & = & -2 \int^{\infty}_{0}
    \left( \tmmathbf{A}_{\alpha} \widehat{\tmmathbf{\nabla}}_{z} P_{t} f,
    \widehat{\tmmathbf{\nabla}}_{z} P_{t} g \right)_{L^{2} ( \hat{T}
    \mathbbm{G} )} \; \mathd  t\\
    & = & -2 \int^{\infty}_{0} \int_{z \in \mathbbm{G}} \Bigg \{ \frac{1}{2}
    \sum_{i=1}^{m} \sum_{\tau = \pm} \alpha_{i}^{x} \; ( X_{i}^{\tau} P_{t} f
    ) ( z )  ( X_{i}^{\tau} P_{t} g ) ( z ) \\
    &  & +  \sum_{j,k=1}^{n} \alpha_{j k}^{y} \; ( Y_{j}
    P_{t} f ) ( z ) ( Y_{k} P_{t} g ) ( z ) \vphantom{\int_{z \in \mathbbm{G}}}   \Bigg \}  \; \mathd
    \mu_{z} ( z ) \; \mathd  t
  \end{eqnarray*}
and the sums and integrals that arise converge absolutely.

\section{Stochastic integrals and martingale transforms}\label{S: stochastic}

In what follows, we assume that we have a complete probability space $( \Omega
, \mathcal{F} ,\mathbbm{P} )$ with a c{\`a}dl{\`a}g (i.e. right continuous
left limit) filtration $( \mathcal{F}_{t} )_{t \geqslant 0}$ of
sub-$\sigma$--algebras of $\mathcal{F}$. We assume as usual that
$\mathcal{F}_{0}$ contains all events of probability zero. All random walks
and martingales are adapted to this filtration.

We define below a continuous-time random process $\mathcal{Z}$ with values in
$\mathbbm{G}$, $\mathcal{Z}_{t} \assign ( \mathcal{X}_{t} , \mathcal{Y}_{t} )
\in \mathbbm{G}_{x} \times \mathbbm{G}_{y}$, having infinitesimal generator
$L= \Delta_{z}$. The pure-jump component $\mathcal{X}_{t}$ is a compound
Poisson jump process on the discrete set $\mathbbm{G}_{x}$, wheras the
continuous component $\mathcal{Y}_{t}$ is a standard brownian motion on the
manifold $\mathbbm{G}_{y}$. Then, It{\^o}'s formula ensures that semi-discrete
``harmonic'' functions $f:\mathbbm{R}^{+} \times \mathbbm{G} \rightarrow
\mathbbm{C}$ solving the backward heat equation $( \partial_{t} + \Delta_{z} )
f=0$ give rise to martingales $M_{t}^{f} \assign f ( t, \mathcal{Z}_{t} )$ for
which we define a class of martingale transforms.

\medskip

\subsection{Stochastic integrals, Martingale transforms and quadratic covariations}\label{SS: stochastic integrals and martingale transforms}

\mbox{}\\

\paragraph{Stochastic integrals on Riemannian manifolds and It{\^o}
integral}Following {\tmname{Emery}} {\cite{Eme2000a}}{\cite{Eme2005a}}, see
also {\tmname{Arcozzi}} {\cite{Arc1995a}}{\cite{Arc1998a}}, we define the
Brownian motion $\mathcal{Y}_{t}$ on $\mathbbm{G}_{y}$, a compact Riemannian
manifold, as the process $\mathcal{Y}_{t} \  : \  \Omega
\rightarrow ( 0,T ) \times \mathbbm{G}_{y}$ such that for all smooth functions
$f:\mathbbm{G}_{y} \rightarrow \mathbbm{R}$, the quantity
\begin{equation}
  f ( \mathcal{Y}_{t} ) -f ( \mathcal{Y}_{0} ) - \frac{1}{2} \int_{0}^{t} (
  \Delta_{y} f ) ( \mathcal{Y}_{s} ) \; \mathd s=: ( I_{\mathd_{y} f} )_{t}
  \label{eq: Ito for continuous processes}
\end{equation}
is an $\mathbbm{R}$--valued continuous martingale. For any adapted continuous
process $\Psi$ with values in the cotangent space $T^{\ast} \mathbbm{G}_{y}$
of $\mathbbm{G}_{y}$, if $\Psi_{t} ( \omega ) \in T^{\ast}_{Y_{t} ( \omega )}
\mathbbm{G}_{y}$ for all $t \geqslant 0$ and $\omega \in \Omega$, then one can
define the {\tmem{continuous}} It{\^o} integral $I_{\Psi}$ of $\Psi$ as
\[ ( I_{\Psi} )_{t} \assign \int_{0}^{t} \langle \Psi_{s} , \mathd
   \mathcal{Y}_{s} \rangle \]
so that in particular
\[ ( I_{\mathd_{y} f} )_{t} \assign \int_{0}^{t} \langle \mathd_{y} f (
   \mathcal{Y}_{s} ) , \mathd \mathcal{Y}_{s} \rangle \]
The integrand above involves the $1$--form of $T_{y}^{\ast} \mathbbm{G}_{y}$
\[ \mathd_{y} f ( y ) \assign \sum_{j} ( Y_{j} f ) ( y ) \ 
   Y_{j}^{\ast} . \]

\paragraph{A pure jump process on $\mathbbm{G}_{x}$}We will now define the
{\tmem{discrete}} $m$--dimensional process $\mathcal{X}_{t}$ on the
{\tmem{discrete}} abelian group $\mathbbm{G}_{x}$ as a generalized compound
Poisson process. In order to do this we need a number of independent variables
and processes:

First, for any given $1 \leqslant i \leqslant m$, let
$\mathcal{N}^{i}_{t}$ be a c{\`a}dl{\`a}g Poisson process of parameter
$\lambda$, that is
\[ \forall t, \hspace{1em} \mathbbm{P} ( \mathcal{N}^{i}_{t} =n ) = \frac{(
   \lambda  t )^{n}}{n!} e^{- \lambda t} . \]
The sequence of instants where the jumps of the $\mathcal{N}^{i}_{t}$ occur is
noted $( T^{i}_{k} )_{k \in \mathbbm{N}}$, with the convention $T^{i}_{0} =0$.

Second, we set
\[ \mathcal{N}_{t} = \sum_{i=1}^{m} \mathcal{N}_{t}^{i} \]
Almost surely, for any two distinct $i$ and $j$, we have $\{ T_{k}^{i} \}_{k
\in \mathbbm{N}} \cap \{ T_{k}^{j} \}_{k \in \mathbbm{N}} = \emptyset$. Let
therefore $\{ T_{k} \}_{k \in \mathbbm{N}} = \cup_{i=1}^{m} \{ T_{k}^{i} \}_{k
\in \mathbbm{N}}$ be the ordered sequence of instants of jumps of
$\mathcal{N}_{t}$ and let $i_{t} \equiv i_{t} ( \omega )$ be the index of the
coordinate where the jump occurs at time $t$. We set $i_{t} =0$ if no jump
occurs. The random variables $i_{t}$ are measurable: $i_{t} = (
\mathcal{N}^{1}_{t} -\mathcal{N}_{t^{-}}^{1} ,\mathcal{N}^{2}_{t}
-\mathcal{N}_{t^{-}}^{2} , \ldots ,\mathcal{N}^{m}_{t}
-\mathcal{N}_{t^{-}}^{m} ) \cdot ( 1,2, \ldots ,m )$. In differential form,
\[ \mathd \mathcal{N}_{t} = \sum_{i=1}^{m} \mathd \mathcal{N}^{i}_{t} = \mathd
   \mathcal{N}^{i_{t}}_{t} . \]

Third, we denote by $( \tau_{k} )_{k \in \mathbbm{N}}$ a sequence of
independent Bernoulli variables
\[ \forall k, \hspace{1em} \mathbbm{P} ( \tau_{k} =1 ) =\mathbbm{P} ( \tau_{k}
   =-1 ) =1/2. \]

Finally, the random walk $\mathcal{X}_{t}$ started at $\mathcal{X}_{0} \in
\mathbbm{G}_{x}$ is the c{\`a}dl{\`a}g compound Poisson process (see e.g.
{\tmname{Protter}} {\cite{Pro2005a}}, {\tmname{Privault}}
{\cite{Pri2009a,Pri2014a}}) defined as
\[ \mathcal{X}_{t} \assign \mathcal{X}_{0} + \sum_{k=1}^{\mathcal{N}_{t}}
   G_{i_{k}}^{\tau_{k}} , \]
where $G_{i}^{\tau} = ( 0, \ldots ,0, \tau g_{i} ,0, \ldots ,0 )$ when $i_{}
\neq 0$ \ and $( 0, \ldots ,0 )$ when $i_{} =0$.

\paragraph{Stochastic integrals on discrete groups}We recall for the
convenience of the reader the derivation of stochastic integrals for jump
processes. We will emphasize the fact that the corresponding It{\^o}'s formula
involves the action of a discrete $1$--form written in a well-chosen local
coordinate system of the discrete {\tmem{augmented}} cotangent plane (see
details below). Let $1 \leqslant k \leqslant \mathcal{N}_{t}$ and let $( T_{k}
,i_{k} , \tau_{k} )$ be respectively the instant, the axis and the direction
of the $k$--th jump. We set $T_{0} =0$. Let $f \assign f ( t,x )$, $t \in
\mathbbm{R}^{+}$, $x \in \mathbbm{G}_{x}$ a function defined on
$\mathbbm{R}^{+} \times \mathbbm{G}_{x}$. Then
\begin{eqnarray*}
  \lefteqn{f ( t, \mathcal{X}_{t} ) - f ( 0,\mathcal{X}_{0} )}\\
  & = & \int_{0}^{t} ( \partial_{t} f ) ( s, \mathcal{X}_{s} ) \mathd s+
  \sum_{i=1}^{m} \int_{0}^{t} ( f ( s, \mathcal{X}_{s} ) -f ( s,
  \mathcal{X}_{s_{-}} ) ) \; \mathd \mathcal{N}^{i}_{s} .
\end{eqnarray*}

At an instant $s$, the integrand in the last term writes as

\begin{eqnarray*}
  \lefteqn{( f ( s, \mathcal{X}_{s} ) -f ( s, \mathcal{X}_{s_{-}} ) )
  \mathd \mathcal{N}^{i}_{s}}\\
   & = & \left( f \left( s, \mathcal{X}_{s_{-}}
  +G_{i}^{\tau_{\mathcal{N}_{s}}} \right) -f ( s, \mathcal{X}_{s_{-}} )
  \right)  \mathd \mathcal{N}^{i}_{s}\\
  & = & \left( X_{i}^{\tau_{\mathcal{N}_{s}}} f \right) ( s,
  \mathcal{X}_{s_{-}} ) \tau_{\mathcal{N}_{s}} \mathd
  \mathcal{N}^{i}_{s}\\
  & = & \frac{1}{2} \left\{ ( X_{i}^{2} f ) ( s, \mathcal{X}_{s_{-}} )
   + \tau_{\mathcal{N}_{s}}  ( X_{i}^{0}
  f ) ( s, \mathcal{X}_{s_{-}} ) \right\} \mathd
  \mathcal{N}^{i}_{s}
\end{eqnarray*}
where we introduced, for all $1 \leqslant i \leqslant m$,
\begin{eqnarray*}
  X_{i}^{0} & \assign & X_{i}^{+} +X_{i}^{-}\\
  X_{i}^{2} & \assign & X_{i}^{+} -X_{i}^{-} .
\end{eqnarray*}
Notice that, for any given $1 \leqslant i \leqslant m$, up to a normalisation
factor, the system of coordinate $( X_{i}^{2} ,X_{i}^{0} )$ is obtained thanks
to a {\tmem{rotation}} of $\pi /4$ of the canonical system of coordinate $(
X_{i}^{+} ,X_{i}^{-} )$. Finally,
\begin{eqnarray*}
 \lefteqn{ f ( t, \mathcal{X}_{t} ) - f ( 0, \mathcal{X}_{0} )}\\
  & = &  \int_{0}^{t} \left\{ ( \partial_{t} f ) ( s, \mathcal{X}_{s} ) +
  \frac{\lambda}{2} ( \Delta_{x} f ) ( s, \mathcal{X}_{s} ) \right\} \mathd s+
  \int_{0}^{t} \left\langle \widehat{\tmmathbf{\mathd}} f ( s,
  \mathcal{X}_{s_{-}} ) , \mathd \widehat{\mathcal{W}}_{s} \right\rangle
  \nonumber\\
  & =: & \int_{0}^{t} \left\{ ( \partial_{t} f ) ( s, \mathcal{X}_{s} ) +
  \frac{\lambda}{2} ( \Delta_{x} f ) ( s, \mathcal{X}_{s} ) \right\} \mathd s+
  \left( I_{\widehat{\tmmathbf{\mathd}}_{x} f} \right)_{t}  . 
\end{eqnarray*}
where we set $ \mathd \mathcal{X}^{i}_{s} =  \tau_{\mathcal{N}_{s}} \mathd \mathcal{N}^{i}_{s}$.
It is easy to see that $  \mathd \mathcal{X} ^{i}_{s}$ is the stochastic
differential of a martingale. Here and in the sequel, we take $\lambda =2$.

\paragraph{Discrete It{\^o} integral}The stochastic integral above shows that
It{\^o} formula (\ref{eq: Ito for continuous processes}) for continuous
processes has a discrete counterpart involving stochastic integrals for jump
processes, namely we have the {\tmem{discrete}} It{\^o} integral
\begin{eqnarray*}
  \left( I_{\widehat{\tmmathbf{\mathd}}_{x} f} \right)_{t} & \assign &
  \frac{1}{2} \sum_{i=1}^{m} \int_{0}^{t} ( X^{2}_{i} f ) ( s,
  \mathcal{X}_{s_{-}} ) \  \mathd ( \mathcal{N}^{i}_{s} - \lambda s
  ) + ( X^{0}_{i} f ) ( s, \mathcal{X}_{s_{-}} ) \  \mathd
  \mathcal{X}_{s}
\end{eqnarray*}
This has a more intrinsic expression similar to the continuous It{\^o}
integral (\ref{eq: Ito for continuous processes}). If we regard the discrete
component $\mathbbm{G}_{x}$ as a ``discrete Riemannian'' manifold, then this
discrete It{\^o} integral involves discrete vectors (resp. $1$--forms) defined
on the {\tmem{augmented}} discrete tangent (resp. cotangent) space
$\hat{T}_{x} \mathbbm{G}_{x}$ (resp. $\hat{T}_{x}^{\ast} \mathbbm{G}_{x}$) of
dimension $2m$ defined as
\begin{eqnarray*}
  \hat{T}_{x} \mathbbm{G}_{x} & = & \tmop{span} \{ X_{1}^{+} ,X_{2}^{+} ,
  \ldots ,X_{1}^{-} ,X_{2}^{-} , \ldots \}\\
  & = & \tmop{span} \{ X_{1}^{2} ,X_{2}^{2} , \ldots ,X_{1}^{0} ,X_{2}^{0} ,
  \ldots \}\\
  \hat{T}_{x}^{\ast} \mathbbm{G}_{x} & = & \tmop{span} \{ ( X_{1}^{+} )^{\ast}
  , ( X_{2}^{+} )^{\ast} , \ldots , ( X_{1}^{-} )^{\ast} , ( X_{2}^{-}
  )^{\ast} , \ldots \}\\
  & = & \tmop{span} \{ ( X_{1}^{2} )^{\ast} , ( X_{2}^{2} )^{\ast} , \ldots ,
  ( X_{1}^{0} )^{\ast} , ( X_{2}^{0} )^{\ast} , \ldots \} .
\end{eqnarray*}
Let $\mathd \widehat{\mathcal{W}}_{s} \in \hat{T}_{\mathcal{X}_{s}}
\mathbbm{G}_{x}$ be the vector and $\widehat{\tmmathbf{\mathd}} f \in
\hat{T}_{\mathcal{X}_{s}}^{\ast} \mathbbm{G}_{x}$ be the $1$--form
respectively defined as:
\begin{equation*}
  \mathd \widehat{\mathcal{W}}_{s} =  \mathd ( \mathcal{N}^{1}_{s} -
  \lambda s )  X_{1}^{2} + \ldots + \mathd ( \mathcal{N}^{m}_{s} -
  \lambda s ) X_{m}^{2} + \mathd
  \mathcal{X}_{s}^{1} X_{1}^{0} + \ldots + \mathd
  \mathcal{X}_{s}^{m}  X_{m}^{0}
\end{equation*}
\begin{equation*}
  \widehat{\tmmathbf{\mathd}}_{x} f  = X_{1}^{2} f  ( X_{1}^{2}
  )^{\ast} + \ldots +X_{m}^{2} f ( X_{m}^{2} )^{\ast} +X_{1}^{0}
  f  ( X_{1}^{0} )^{\ast} + \ldots +X_{m}^{0} f  (
  X_{m}^{0} )^{\ast}
\end{equation*}
We have with these notations
\begin{eqnarray*}
  \left( I_{\widehat{\tmmathbf{\mathd^{}}}_{x} f} \right)_{t} & \assign &
  \left\langle \widehat{\tmmathbf{\mathd}}_{x} f, \mathd
  \widehat{\mathcal{W}}_{s} \right\rangle_{\hat{T}_{x}^{\ast} \mathbbm{G}_{x}
  \times \hat{T}_{x} \mathbbm{G}_{x}}
\end{eqnarray*}
where the factor $1/2$ is included in the pairing $\langle \cdot , \cdot
\rangle_{\hat{T}_{x}^{\ast} \mathbbm{G}_{x} \times \hat{T}_{x}
\mathbbm{G}_{x}}$.

\paragraph{Semi--discrete stochastic integrals}Let finally $\mathcal{Z}_{t} = (
\mathcal{X}_{t} , \mathcal{Y}_{t} )$ be a semi-discrete random walk on the
cartesian product $\mathbbm{G}=\mathbbm{G}_{x} \times \mathbbm{G}_{y}$, where
$\mathcal{X}_{t}$ is the random walk above defined on $\mathbbm{G}_{x}$ with
generator $\Delta_{x}$ and where $\mathcal{Y}_{t}$ is the Brownian motion
defined on $\mathbbm{G}_{y}$ with generator $\Delta_{y}$. For $f \assign f (
t,z ) =f ( t,x,y )$ defined from $\mathbbm{R}^{+} \times \mathbbm{G}$ onto
$\mathbbm{C}$, we have easily the stochastic integral involving both discrete
and continuous parts:
\begin{eqnarray*}
  f ( t, \mathcal{Z}_{t} ) & = & \int_{0}^{t} \{ ( \partial_{t} f ) (
  s,\mathcal{Z}_{s} ) + ( \Delta_{z} f ) ( s,\mathcal{Z}_{s} ) \} \; \mathd s+
  \left( I_{\widehat{\tmmathbf{\mathd}}_{z} f} \right)_{t}
\end{eqnarray*}
where the {\tmem{semi-discrete}} It{\^o} integral writes as
\begin{eqnarray*}
  \left( I_{\widehat{\tmmathbf{\mathd}}_{z} f} \right)_{t} & \assign & \left(
  I_{\widehat{\tmmathbf{\mathd}}_{x} f} \right)_{t} + \left(
  I_{\tmmathbf{\mathd}_{y} f} \right)_{t}\\
  & \assign & \int_{0}^{t} \left\langle \widehat{\tmmathbf{\mathd}}_{x} f (
  s, \mathcal{Z}_{s_{-}} ) , \mathd \widehat{\mathcal{W}}_{s}
  \right\rangle_{\hat{T}_{\mathcal{X}_{s}}^{\ast} \mathbbm{G}_{x} \times
  \hat{T}_{\mathcal{X}_{s}} \mathbbm{G}_{x}} \\
  &&+ \int_{0}^{t} \left\langle
  \tmmathbf{\mathd}_{y} f ( s, \mathcal{Z}_{s_{-}} ) , \mathd \mathcal{Y}_{s}
  \right\rangle_{\hat{T}_{\mathcal{Y}_{s}}^{\ast} \mathbbm{G}_{y} \times
  \hat{T}_{\mathcal{Y}_{s}} \mathbbm{G}_{y}}.
\end{eqnarray*}

\mbox{}\\

\paragraph{Martingale transforms} We are interested in martingale transforms
allowing us to represent second order Riesz transforms. Let $f ( t,z )$ be a
solution to the heat equation $\partial_{t} - \Delta_{z} =0$. Fix $T>0$ and
$\mathcal{Z}_{0} \in \mathbbm{G}$. Then define
\[ \forall 0 \leqslant t \leqslant T, \hspace{1em} M_{t}^{f,T,\mathcal{Z}_{0}}
   =f ( T-t,\mathcal{Z}_{t} ) . \]
This is a martingale since $f ( T-t )$ solves the backward heat equation
$\partial_{t} + \Delta_{z} =0$, and we have in terms of stochastic integrals
\[ M_{t}^{f,T,\mathcal{Z}_{0}} =f ( T-t,\mathcal{Z}_{t} ) =f (
   T,\mathcal{Z}_{0} ) + \int_{0}^{t} \left\langle
   \widehat{\tmmathbf{\mathd}}_{z} f ( T-s,\mathcal{Z}_{s_{-}} ) , \mathd
   \mathcal{Z}_{s} \right\rangle \]
Given $\tmmathbf{A}_{\alpha}$ the $\mathbbm{C}^{( 2m+n ) \times ( 2m+n )}$
matrix defined earlier, we note $M_{t}^{\alpha ,f,T,\mathcal{Z}_{0}}$ the
martingale transform $\tmmathbf{A}_{\alpha} \ast M_{t}^{f,T,\mathcal{Z}_{0}}$
defined as
\begin{eqnarray*}
  M_{t}^{\alpha ,f,T,\mathcal{Z}_{0}} & \assign & f ( T,\mathcal{Z}_{0} ) +
  \int_{0}^{t} \left( \tmmathbf{A}_{\alpha} \widehat{\tmmathbf{\nabla}}_{z} f
  ( s,\mathcal{Z}_{s_{-}} ) , \mathd \mathcal{Z}_{s} \right)\\
  & = & f ( T,\mathcal{Z}_{0} ) + \int_{0}^{t} \left\langle
  \widehat{\tmmathbf{\mathd}}_{z} f ( T-s,\mathcal{Z}_{s_{-}} )
  \tmmathbf{A}_{\alpha}^{\ast} , \mathd \mathcal{Z}_{s} \right\rangle
\end{eqnarray*}
where the first integral involves the $L^{2}$ scalar product on $\hat{T}_{z}
\mathbbm{G} \times \hat{T}_{z} \mathbbm{G}$ and the second integral involves
the duality $\hat{T}_{z}^{\ast} \mathbbm{G} \times \hat{T}_{z} \mathbbm{G}$.
In differential form:
\begin{eqnarray*}
  \lefteqn{\mathd M_{t}^{\alpha ,f,T,\mathcal{Z}_{0}} }\\
  & = & \left(
  \tmmathbf{A}_{\alpha} \widehat{\tmmathbf{\nabla}}_{z} f (
  s,\mathcal{Z}_{s_{-}} ) , \mathd \mathcal{Z}_{s} \right)\\
  & = & \sum_{i=1}^{m} \alpha_{i}^{x}  \left\{ ( X^{2}_{i} f ) (
  T-t,\mathcal{Z}_{t_{-}} ) \; \mathd ( \mathcal{N}^{i}_{t} - \lambda t ) + (
  X^{0}_{i} f ) ( t,\mathcal{Z}_{t_{-}} ) \; \mathd \mathcal{X}_{t}^{i}
  \right\}\\
  &  & + \sum_{j=1}^{n} \alpha_{j,k}^{y} \; ( X_{j} f ) (
  T-t,\mathcal{Z}_{t_{-}} ) \; \mathd \mathcal{Y}_{t}^{k}
\end{eqnarray*}

\paragraph{Quadratic covariation and subordination}We have the quadratic
covariations (see {\tmname{Protter }}{\cite{Pro2005a}},
{\tmname{Dellacherie}}--{\tmname{Meyer}} {\cite{DelMey1982a}}, or
{\tmname{Privault}} {\cite{Pri2009a,Pri2014a}}). Since
\begin{eqnarray*}
  \mathd [ \mathcal{N}^{i} - \lambda t, \mathcal{N}^{i} - \lambda t ]_{t} & =
  & \mathd \mathcal{N}_{t}^{i}\\
  \mathd [ \mathcal{N}^{i} - \lambda t, \mathcal{X}^{i} ]_{t} & = &
  \tau_{\mathcal{N}_{t}} \  \mathd \mathcal{N}_{t}^{i}\\
  \mathd [ \mathcal{X}^{i} , \mathcal{X}^{i} ]_{t} & = & \mathd
  \mathcal{N}_{t}^{i}\\
  \mathd [ \mathcal{Y}^{j} , \mathcal{Y}^{j} ]_{t} & = & \mathd t,
\end{eqnarray*}

it follows that

\begin{eqnarray} \label{eq:quadratic_covariation_of_f_and_g}
  \mathd [ M^{f} ,M^{g} ]_{t} &=& \sum_{i=1}^{m} \sum_{\tau = \pm} 
  ( X^{\tau}_{i} f ) \; ( X^{\tau}_{i} g ) ( T-t, \mathcal{Z}_{t_{-}} )
   \mathbbm{1} ( \tau_{\mathcal{N}_{t}} = \tau 1 )   \mathd
  \mathcal{N}^{i}_{t} \\
  &&+ \left( \tmmathbf{\nabla}_{y} f, \tmmathbf{\nabla}_{y}
  g \right) ( T-t, \mathcal{Z}_{t_{-}} )  \mathd t \nonumber
\end{eqnarray}

\paragraph{Differential subordination}Following {\tmname{Wang}}
{\cite{Wan1995a}}, given two adapted c{\`a}dl{\`a}g Hilbert space valued
martingales $X_{t}$ and $Y_{t}$, we say that $Y_{t}$ is differentially
subordinate by quadratic variation to $X_{t}$ if $| Y_{0} |_{\mathbbm{H}}
\leqslant | X_{0} |_{\mathbbm{H}}$ and $[ Y,Y ]_{t} - [ X,X ]_{t}$ is
nondecreasing nonnegative for all $t$. In our case, we have
\begin{eqnarray*}
  \mathd [ M^{\alpha ,f} ,M^{\alpha ,f} ]_{t} & = & \sum_{i=1}^{m} |
  \alpha_{i}^{x} |^{2} \; \left\{ ( X^{+}_{i} f )^{2} ( T-t,
  \mathcal{Z}_{t_{-}} )  \mathbbm{1} ( \tau_{\mathcal{N}_{t}} =1 )
  \right.\\
  &  & +  \left. ( X^{-}_{i} f )^{2} ( T-t,
  \mathcal{Z}_{t_{-}} )  \mathbbm{1} ( \tau_{\mathcal{N}_{t}} =-1
  ) \right\}  \mathd \mathcal{N}^{i}_{t}\\
  &  &  + \left( \tmmathbf{A}^{y}_{\alpha}  
  \tmmathbf{\nabla}_{y} f,\tmmathbf{A}^{y}_{\alpha}   \tmmathbf{\nabla}_{y} f
  \right) ( T-t, \mathcal{Z}_{t_{-}} ) \mathd t.
\end{eqnarray*}
Hence
\begin{equation}
  \mathd [ M^{\alpha ,f} ,M^{\alpha ,f} ]_{t} \leqslant \|
  \tmmathbf{A}_{\alpha} \|^{2}_{2}  \  \mathd [ M^{f} ,M^{f} ]_{t}
  \label{eq: differential subordination} .
\end{equation}
This means that $M^{\alpha ,f}_{t}$ is differentially subordinate to $\|
\tmmathbf{A}_{\alpha} \|_{2} M_{t}^{f}$.

\subsection{Martingale inequalities under differential subordination}\label{SS: martingale inequalities} In the final part of the section we discuss a number of sharp martingale inequalities which hold under the assumption of the differential subordination imposed on the processes. Our starting point is the following celebrated $L^p$ bound.

\begin{theorem}\label{T: Wang}\dueto{Wang, 1995}
Suppose that $X$ and $Y$ are martingales taking values in a Hilbert space $\mathbbm{H}$ such that $Y$ is differentially subordinate to $X$. Then for any $1<p<\infty$ we have
$$ ||Y||_p\leq (p^\ast-1)||X||_p$$
and the constant $p^\ast-1$ is the best possible, even if $\mathbbm{H}=\mathbbm{R}$.
\end{theorem}

This result was first proved by {\tmname{Burkholder}} in \cite{Bur1984a} in the following discrete-time setting. Suppose that $(X_n)_{n\geq 0}$ is an $\mathbbm{H}$-valued martingale and $(\alpha_n)_{n\geq 0}$ is a predictable sequence with values in $[-1,1]$. Let $Y:=\alpha*X$ be the martingale transform of $X$ defined for almost all $\omega\in \Omega$ by
$$ Y_0(\omega)=\alpha_0 X_0(\omega)\quad \mbox{ and }\quad (Y_{n+1}-Y_n)(\omega)=\alpha_n(X_{n+1}-X_n)(\omega).$$
Then the above $L^p$ bound holds true and the constant $p^\ast-1$ is optimal. The general continuous-time version formulated above is due to {\tmname{Wang}} \cite{Wan1995a}. To see that the preceding discrete-time version is indeed a special case, treat a discrete-time martingale $(X_n)_{n\geq 0}$ and its transform $(Y_n)_{n\geq 0}$ as continuous-time processes via $X_t=X_{\lfloor t\rfloor}$, $Y_t=Y_{\lfloor t\rfloor}$ for $t\geq 0$; then $Y$ is differentially subordinate to $X$.

In 1992, {\tmname{Choi}} \cite{Cho1992a} established the following non-symmetric, discrete-time version of the $L^p$ estimate.

\begin{theorem}
 {\dueto{Choi, 1992}}Suppose that $(X_n)_{n\geq 0}$ is a real-valued discrete time martingale and let $(Y_n)_{n\geq 0}$ be its transform by a predictable sequence $(\alpha_n)_{n\geq 0}$ taking values in $[0,1]$. 
  Then there exists a constant $\mathfrak{C}_{p}$ depending only on $p$ such
  that $\| Y \|_{p} \leqslant \mathfrak{C}_{p} \| X \|_{p}$ and the estimate
  is best possible.
\end{theorem}

This result can be regarded as a non-symmetric version of the previous theorem, since the transforming sequence $(\alpha_n)_{n\geq 0}$ takes values in a non-symmetric interval $[0,1]$. 
There is a natural question whether the estimate can be extended to the continuous-time setting; in particular, this gives rise to the problem of defining an appropriate notion of non-symmetric differential subordination. The following statement obtained by 
{\tmname{Ba{\~n}uelos}} and {\tmname{Os\c ekowski}}  addresses both these questions. For any real numbers $a<b$ and any $1<p<\infty$, let $\mathfrak{C}_{a,b,p}$ be the constant introduced in {\cite{BanOse2012a}}.  

\begin{theorem}
  \label{L: Choi Lp norms of martingales}{\dueto{Banuelos--Os\c ekowski,
  2012}}Let $(X_{t})_{t\geq 0}$ and $(Y_{t})_{t\geq 0}$ be two real-valued martingales satisfying
  \begin{equation}\label{nonsymmetric_ds}
 \mathd \left[ Y- \frac{a+b}{2} X,Y- \frac{a+b}{2} X \right]_{t} \leqslant
     \mathd \left[ \frac{b-a}{2} X, \frac{b-a}{2} X \right]_{t} 
     \end{equation}
  for all $t \geqslant 0$. Then for all $1<p< \infty$, we have $\| Y \|_{p}
  \leqslant \mathfrak{C}_{a,b,p} \| X \|_{p}$.
\end{theorem}
The condition \eqref{nonsymmetric_ds} is the continuous counterpart of the condition that the transforming sequence $(\alpha_n)_{n\geq 0}$ takes values in the interval $[a,b]$. Thus, in particular, Choi's constant $\mathfrak{C}_p$ is, in the terminology of the above theorem, equal to $\mathfrak{C}_{p,0,1}$.

We return to the context of the ``classical'' differential subordination introduced in the preceding subsection and study other types of martingale inequalities.
The following statements, obtained by {\tmname{Ba\~nuelos}--\tmname{Os\c ekowski}}, \cite{BanOse2015a} will allow us to deduce sharp weak-type and logarithmic estimates for Riesz transforms, respectively. 

\begin{theorem}\label{T: weak-type martingales}{\dueto{Banuelos--Os\c ekowski,
  2015}}
Suppose that $X$ and $Y$ are martingales taking values in a Hilbert space $\mathbbm{H}$ such that $Y$ is differentially subordinate to $X$. 

(i) Let $1<p<2$. Then for any $t\geq 0$,
$$ \mathbb{E}\max\left\{|Y_t|-\frac{p^{-1/(p-1)}}{2}\Gamma\left(\frac{p}{p-1}\right),0\right\}\leq \mathbb{E}|X_t|^p.$$

(ii) Suppose that $2<p<\infty$. Then for any $t\geq 0$,
$$ \mathbb{E}\max\left\{|Y_t|-1+p^{-1},0\right\}\leq \frac{p^{p-2}}{2}\mathbb{E}|X_t|^p.$$

\noindent Both estimates are sharp: for each $p$, the numbers $\frac{p^{-1/(p-1)}}{2}\Gamma\left(\frac{p}{p-1}\right)$ and $1-p^{-1}$ cannot be decreased.
\end{theorem}

Recall that $\Phi,\,\Psi:[0,\infty)\to [0,\infty)$ are conjugate Young functions given by $\Phi(t)=e^t-1-t$ and $\Psi(t)=(t+1)\log(t+1)-t$.

\begin{theorem}{\dueto{Banuelos--Os\c ekowski,
  2015}}
Suppose that $X$ and $Y$ are martingales taking values in a Hilbert space $\mathbbm{H}$ such that $Y$ is differentially subordinate to $X$. Then for any $K>1$ and any $t\geq 0$ we have
$$ \mathbb{E}\left\{|Y_t|-(2(K-1))^{-1},0\right\}\leq K \mathbb{E}\Psi(|X_t|).$$
For each $K$, the constant $(2(K-1))^{-1}$ appearing on the left, is the best possible (it cannot be replaced by any smaller number).
\end{theorem}

The following exponential estimate, established by {\tmname{Os\c ekowski}} in \cite{Ose2013a}, can be regarded as a dual statement to the above logarithmic bound. 

\begin{theorem}{\dueto{Os\c ekowski,
  2013}}
Assume that $X$, $Y$ are $\mathbbm{H}$-valued martingales such that $||X||_\infty\leq 1$ and $Y$ is differentially subordinate to $X$. Then for any $K>1$ and any $t\geq 0$ we have
\begin{equation}
\mathbb{E} \Phi(|Y_t|/K)\leq \frac{1}{2K(K-1)}\mathbb{E} |X_t|.
\end{equation}
\end{theorem}

Finally, we will need the following sharp $L^q\to L^p$ estimate, established by {\tmname{Os\c ekowski}} in \cite{Ose2014b}, which will allow us to deduce the corresponding estimate for Riesz transforms.

\begin{theorem}{\dueto{Os\c ekowski,
  2014}}
Assume that $X$, $Y$ are $\mathbbm{H}$-valued martingales such that $Y$ is differentially subordinate to $X$. Then for any $1\leq p<q<\infty$ there is a constant $L_{p,q}$ such that 
\begin{equation}
 \mathbb{E} \max\{|Y_t|^p-L_{p,q},0\}\leq \mathbb{E}|X_t|^q.
\end{equation}
\end{theorem}
Actually, the paper \cite{Ose2014b} identifies, for any $p$ and $q$ as above, the optimal (i.e., the least) value of the constant $L_{p,q}$ in the estimate above. As the description of this constant is a little complicated (and will not be needed in our considerations below), we refer the reader to that paper for the formal definition of $L_{p,q}$.

Let us conclude with the observation which will be crucial in the proofs of our main results. Namely, all the martingale inequalities presented above are of the form $\mathbb{E}\zeta(|Y_t|)\leq \mathbb{E}\xi(|X_t|)$, $t\geq 0$, where $\zeta$, $\xi$ are certain convex functions. This will allow us to successfully apply a conditional version of Jensen's inequality.

\section{Proofs of the main results}\label{S: proofs of the main results}

We turn our attention to the proofs of the estimates for $R_\alpha^2$ formulated in the introductory section. We will focus on Theorems \ref{T: p minus 1 estimate}, \ref{T: Choi constant estimate} and \ref{T: weak-type estimate} only; the remaining statements are established by similar arguments. Also, we postpone the proof of the sharpness of these estimates to the next section.

\subsection{Proof of Theorem \ref{T: p minus 1 estimate}}Recall that the
subordination estimate (\ref{eq: differential subordination}) shows that the
martingale transform $Y_{t} \assign M_{t}^{\alpha}$ is differentially
subordinate to the martingale $X_{t} \assign \| \tmmathbf{A}_{\alpha} \|_{2}
\hspace{.1em} M_{t}^{f}$. Therefore, by Theorem \ref{T: Wang}, 
we immediately obtain that
\[ \| M_{t}^{\alpha ,f} \|_{p} \leqslant \|
     \tmmathbf{A}_{\alpha} \|_{2} \hspace{.1em} ( p^{\ast} -1 ) \hspace{.1em} \|
     M_{t}^{f} \|_{p}  \]
for all $t\geq 0$. Since the operator $\mathcal{T}^{\alpha}$ is a conditional expectation
of $M_{t}^{\alpha ,f}$, an application of Jensen's inequality proves the estimate $\| \mathcal{T}^{\alpha}
\|_{p} \leqslant \| \tmmathbf{A}_{\alpha} \|_{2} \hspace{.1em} ( p^{\ast} -1
)$, which is the desired bound.

\subsection{Proof of Theorem \ref{T: Choi constant estimate}}The argument is the same as above and exploits the fine-tuned $L^p$ estimate of Theorem \ref{L: Choi Lp norms of martingales} applied to $X_{t} =M_{t}^{f}$ and $Y_{t} =M_{t}^{\alpha ,f}$. It is not
difficult to prove that the difference of quadratic variations above writes in
terms of a jump part and a continuous part as
 \begin{eqnarray*}
     \lefteqn{ \left[ Y- \frac{a+b}{2} X,Y- \frac{a+b}{2} X \right]_{t} - \mathd
     \left[ \frac{b-a}{2} X, \frac{b-a}{2} X \right]_{t} }\\
     &=& \sum_{i=1}^{m} \sum_{\pm} ( \alpha_{i}^{x} -a ) (
     \alpha_{i}^{x} -b ) \hspace{.1em} ( X^{\pm}_{i} f )^{2} ( \mathcal{B}_{t}
     ) \hspace{.1em} \mathbbm{1} ( \tau_{N_{t}} = \pm 1 ) \hspace{.1em} \mathd
     \mathcal{N}_{t}^{i}\\
     &&+ \left\langle \left( \tmmathbf{A}_{\alpha}^{y} -a
     \tmmathbf{I} \right) \left( \tmmathbf{A}_{\alpha}^{y} -b \tmmathbf{I}
     \right) \hspace{.1em} \tmmathbf{\nabla}_{y} f ( \mathcal{B}_{t} )
     ,  \tmmathbf{\nabla}_{y} f ( \mathcal{B}_{t} )
     \right\rangle \hspace{.1em} \mathd t,
   \end{eqnarray*} 
which is nonpositive since we assumed precisely $a \tmmathbf{I} \leqslant
\tmmathbf{A}_{\alpha} \leqslant b \tmmathbf{I}$. Thus, the estimate of Theorem \ref{T:
Choi constant estimate} follows. The sharpness is established in a similar manner.{\hfill}$\Box$

\subsection{Proof of Theorem \ref{T: weak-type estimate}} We will focus on the case $1<p<2$; for remaining values of $p$ the argument is similar. An application of Theorem \ref{T: weak-type martingales} to the processes $X_{t} = \| \tmmathbf{A}_{\alpha} \|_{2}
\hspace{.1em}M_{t}^{f}$ and $Y_{t} =M_{t}^{\alpha ,f}$ yields
$$ \mathbb{E}\max\left\{|M_t^{\alpha,f}|-\frac{p^{-1/(p-1)}}{2}\Gamma\left(\frac{p}{p-1}\right),0\right\}\leq  \| \tmmathbf{A}_{\alpha} \|_{2}^p
\hspace{.1em}\mathbb{E}|M_t^f|^p$$
and hence, by Jensen's inequality, we obtain
$$ \int_{\mathbbm{G}}\max\left\{|R_\alpha^2f|-\frac{p^{-1/(p-1)}}{2}\Gamma\left(\frac{p}{p-1}\right),0\right\}d\mu_z\leq  \| \tmmathbf{A}_{\alpha} \|^p_{2}
\hspace{.1em}||f||_{L^p(\mathbbm{G})}^p.$$
Therefore, if $E$ is an arbitrary measurable subset of $\mathbbm{G}$, we get
\begin{align*}
 \int_E |R_\alpha^2f|d\mu_z&\leq \int_E \left(|R_\alpha^2f|-\frac{p^{-1/(p-1)}}{2}\Gamma\left(\frac{p}{p-1}\right)\right)\mbox{d}\mu_z\\
&\quad +\frac{p^{-1/(p-1)}}{2}\Gamma\left(\frac{p}{p-1}\right)\mu_z(E)\\
&\leq ||f||_{L^p(\mathbbm{G})}^p+\frac{p^{-1/(p-1)}}{2}\Gamma\left(\frac{p}{p-1}\right)\mu_z(E).
\end{align*}
Apply this bound to $\lambda f$, where $\lambda$ is a nonnegative parameter, then divide both sides by $\lambda$ and optimize the right-hand side over $\lambda$ to get the desired assertion.

\section{Sharpness}

The proof of the sharpness of the different results is made in several steps. In some cases the sharpness for certain
second order Riesz transform estimates in the continuous setting (such as in Theorem \ref{T: LpRieszsemidiscrete}) is already known. In these cases we prove below the sharpness for the discrete (or semidiscrete) case by using sequences of finite difference approximates of continuous functions and their finite difference second order Riesz transforms. In other cases, we need to prove first sharpness for certain continuous second order Riesz transforms. The key point here is to transfer the sharp result for zigzag martingales into a sharp result for certain continuous second order Riesz transforms by the laminate technique. We will illustrate this for the weak-type estimate of Theorem \ref{T: weak-type estimate} and establish  the following statement.

\begin{theorem}\label{laminate}
Let $\Theta:[0,\infty)\to [0,\infty)$ be a given function and let $\lambda>0$ be a fixed number. Assume further that there is a pair $(F,G)$ of finite martingales starting from $(0,0)$ such that $G$ is a $\pm 1$-transform of $F$ and
$$
 \E (|G_\infty|-\lambda)_+>\E \Theta(|F_\infty|).
$$
Then there is a function $f:\R^2\to \R$ supported on the unit disc $\mathbb{D}$ of $\R^2$ such that
$$
 \int_{\R^2}\big(|(R_1^2-R_2^2)f|-\lambda\big)_+\mbox{d}x>\int_{\mathbb{D}} \Theta(|f|)\mbox{d}x.
$$
\end{theorem}

We will prove this statement with the use of laminates, important family of probability measures on matrices. It is convenient to split this section into several separate parts.  For the sake of convenience, and to make this section as self contained as possible, we recall the preliminaries on laminates and their connections to martingales from \cite{BSV} and \cite{O5}, Section 4.2. 

\subsection{Laminates} Assume that $\R^{m\times n}$ stands for the space of all real matrices of dimension  $m\times n$ and  $\R^{n\times n}_{sym}$ denote  the subclass of $\R^{n\times n}$ which consists of all symmetric matrices of dimension $n\times n$.

\begin{dfn}
A function $f:\R^{m\times n} \to \R$ is said to be {\it rank-one convex}, if  for all $A,B \in \R^{m\times n}$ with $\textrm{rank }B= 1$, the function 
$t\mapsto f(A+tB)$ is convex
\end{dfn}

For other equivalent definitions of rank-one convexity, see \cite[p.~100]{Dac}.  Suppose that  $\mathcal{P}=\mathcal{P}(\R^{m\times n})$ is the class of all compactly supported probability measures on $\R^{m \times n}$.
For a measure $\nu \in \mathcal{P}$, we define $$\overline{\nu} = \int_{\R^{m\times n}}{X d\nu(X)},$$
the associated \emph{center of mass} or \textit{barycenter} of $\nu.$

\begin{dfn}
We say that a measure $\nu \in \mathcal{P}$ is a \textit{laminate} (and write $\nu\in\mathcal{L}$), if 
\begin{equation*}
f(\overline{\nu}) \leq \int_{\R^{m\times n}}f \mbox{d}\nu
\end{equation*} 
for all rank-one convex functions $f$. The set of laminates with barycenter $0$ is denoted by $\mathcal{L}_0(\R^{m\times n})$. 
\end{dfn}

Laminates can be used to obtain lower bounds for solutions of certain PDEs, as observed by  Faraco in \cite{F}. In addition, laminates appear naturally in the context of convex integration, where they lead to interesting counterexamples, see e.g. \cite{AFS}, \cite{CFM}, \cite{KMS}, \cite{MS99} and \cite{SzCI}. For our results here we will be interested in the case of $2\times 2$ symmetric matrices. The key observation is  that laminates can be regarded as probability measures that record the distribution of the gradients of smooth maps: see Corollary \ref{coro} below. We briefly explain this and refer the reader to the works \cite{Kirchheim}, \cite{MS99} and \cite{SzCI} for full details.  

\begin{dfn}
Let $U$ be a subset of $\R^{2\times 2}$ and let $\mathcal{PL}(U)$ denote  the smallest
class of probability measures on $U$ which 

\begin{itemize}

\item[(i)] contains all measures of the form $\lambda \delta_A+(1-\lambda)\delta_B$ with $\lambda\in [0,1]$ and satisfying $\textrm{rank}(A-B)=1$;

\item[(ii)] is closed under splitting in the following sense: if $\lambda\delta_A+(1-\lambda)\nu$ belongs to $\mathcal{PL}(U)$ for some $\nu\in\mathcal{P}(\R^{2\times 2})$ and $\mu$ also belongs to $\mathcal{PL}(U)$ with $\overline{\mu}=A$, then also $\lambda\mu+(1-\lambda)\nu$ belongs to $\mathcal{PL}(U)$.
\end{itemize}

The class $\mathcal{PL}(U)$ is called the  \emph{prelaminates} in $U$. 
\end{dfn} 

It follows immediately from the definition that the class $\mathcal{PL}(U)$ only contains atomic measures. Also, by a successive application of Jensen's inequality, we have the inclusion $\mathcal{PL}\subset\mathcal{L}$. The following are two well known lemmas in the theory of laminates; see \cite{AFS}, \cite{Kirchheim}, \cite{MS99}, \cite{SzCI}. 

\begin{lemma}
Let $\nu=\sum_{i=1}^N\lambda_i\delta_{A_i}\in\mathcal{PL}(\R^{2\times 2}_{sym})$ with $\overline{\nu}=0$. Moreover, let
$0<r<\tfrac{1}{2}\min|A_i-A_j|$ and $\delta>0$. For any bounded domain $\mathcal{B}\subset\R^2$ there exists $u\in W^{2,\infty}_0(\mathcal{B})$ such that $\|u\|_{C^1}<\delta$ and for all $i=1\dots N$
$$
\left\| \{x\in\mathcal{B}:\,|D^2u(x)-A_i|<r\} \right\| = \lambda_i |\mathcal{B}|.
$$
\end{lemma}

\begin{lemma}
Let $K\subset\R^{2\times 2}_{sym}$ be a compact convex set and suppose that $\nu\in\mathcal{L}(\R^{2\times 2}_{sym})$ satisfies $\operatorname*{supp}\nu\subset K$. For any relatively open set $U\subset\R^{2\times 2}_{sym}$ 
with $K\subset U$,  there exists a sequence $\nu_j\in \mathcal{PL}(U)$ of prelaminates with $\overline{\nu}_j=\overline{\nu}$ and $\nu_j\overset{*}{\rightharpoonup}\nu$, where $\overset{*}{\rightharpoonup}$ denotes weak convergence of measures. 
\end{lemma}

Combining these two lemmas and using a simple mollification, we obtain the following statement, proved by Boros, Sh\'ekelyhidi Jr. and Volberg \cite{BSV}. It exhibits the connection  between laminates supported on symmetric matrices and second derivatives of functions.  It will be our main tool in the proof of the sharpness. Recall that $\mathbb{D}$ denotes the unit disc of $\mathbb{C}$.

\begin{cor}\label{coro}
Let $\nu\in\mathcal{L}_0(\R^{2\times 2}_{sym})$. Then there exists a sequence $u_j\in C_0^{\infty}(\mathbb{D})$ with uniformly bounded second derivatives, such that
$$
\frac{1}{|\mathbb{D}|}\int_{\mathbb{D}} \phi(D^2u_j(x))\,\mbox{d}x\,\to\,\int_{\R^{2\times 2}_{sym}}\phi\,\mbox{d}\nu
$$
for all continuous $\phi:\R^{2\times 2}_{sym}\to\R$. 
\end{cor}

\subsection{Biconvex functions and a special laminate} 
The next step in our analysis is devoted to the introduction of a certain special laminate. We need some additional notation. 
A function $\zeta:\R\times \R\to \R$ is said to be \emph{biconvex} if for any fixed $z\in \R$, the functions $x\mapsto \zeta(x,z)$ and $y\mapsto \zeta(z,y)$ are convex. Now, take the martingales $F$ and $G$ appearing in the statement of Theorem \ref{laminate}. Then the martingale pair
$$ (\mathtt{F},\mathtt{G}):=\left(\frac{F+G}{2},\frac{F-G}{2}\right)$$
is finite, starts from $(0,0)$ and has the following \emph{zigzag} property: for any $n\geq 0$ we have $\mathtt{F}_n=\mathtt{F}_{n+1}$ with probability $1$ or $\mathtt{G}_n=\mathtt{G}_{n+1}$ almost surely; that is, in each step $(\mathtt{F},\mathtt{G})$ moves either vertically, or horizontally. Indeed, this follows directly from the assumption that $G$ is a $\pm 1$-transform of $F$. This property combines nicely with biconvex functions: if $\zeta$ is such a function, then a successive application of Jensen's inequality gives
\begin{equation}\label{dm} 
\E \zeta(\mathtt{F}_{n},\mathtt{G}_{n})\geq \E\zeta(\mathtt{F}_{n-1},\mathtt{G}_{n-1})\geq \ldots\geq \E \zeta(\mathtt{F}_0,\mathtt{G}_0)=\zeta(0,0).
\end{equation}

The distribution of the terminal variable $(\mathtt{F}_{\infty},\mathtt{G}_{\infty})$ gives rise to a probability measure $\nu$ on $\R^{2\times 2}_{sym}$: put
$$ \nu\left(\diag(x,y)\right)=\mathbb{P}\big((\mathtt{F}_{\infty},\mathtt{G}_{\infty})=(x,y)\big),\qquad (x,y)\in \R^2,$$
where  $\diag(x,y)$ stands for the diagonal matrix 
$ \left(\begin{array}{cc}
x & 0\\ 
0 & y
\end{array}\right).$ 
Observe  that $\nu$ is a laminate of barycenter $0$. Indeed, if $\psi:\R^{2\times 2}\to \R$ is a rank-one convex, then $(x,y)\mapsto \psi(\operatorname*{diag}(x,y))$ is biconvex and thus, by \eqref{dm}, 
\begin{align*}
 \int_{\R^{2\times 2}}\psi\mbox{d}\nu
&=\E \psi(\operatorname*{diag}(\mathtt{F}_{\infty},\mathtt{G}_{\infty})) \geq \psi(\operatorname*{diag}(0,0))=\psi(\bar{\nu}).
\end{align*}
Here we used the fact that $(\mathtt{F},\mathtt{G})$ is finite, so $(\mathtt{F}_\infty,\mathtt{G}_\infty)=(\mathtt{F}_n,\mathtt{G}_n)$ for some $n$. 

\subsection{A proof of Theorem \ref{laminate}} Consider a continuous function $\phi:\R^{2\times 2}_{sym}\to \R$ given by 
$$\phi(A)=(|A_{11}-A_{22}|-\lambda)_+-\Theta(|A_{11}+A_{22}|).$$
 By Corollary \ref{coro}, there is a functional sequence $(u_j)_{j\geq 1}\subset C_0^\infty(\mathbb{D})$ such that
\begin{align*}
 \frac{1}{|\mathbb{D}|}\int_{\R^2}\phi(D^2u_j)\mbox{d}x&=\frac{1}{|\mathbb{D}|}\int_{\mathbb{D}}\phi(D^2u_j)\mbox{d}x\\
&\xrightarrow{j\to\infty}\int_{\R^{2\times 2}_{sym}}\phi\mbox{d}\nu
=\E(|G_\infty|-\lambda)_+-\E \Theta(|F_\infty|)>0.
\end{align*}
Therefore, for sufficiently large $j$, we have
$$ \int_{\R^2}\left(\left\| \frac{\partial^2 u_j}{\partial x^2}-\frac{\partial^2 u_j}{\partial y^2}\right\|-\lambda\right)_+\mbox{d}x\mbox{d}y>\int_{\R^2}\Theta\left(|\Delta u_j|\right)\mbox{d}x\mbox{d}y.$$  
 Setting $f=\Delta u_j$, we obtain the desired assertion. 

\smallskip

In the remaining part of this subsection, let us briefly explain how Theorem \ref{laminate} yields the sharpness of weak-type and logarithmic estimates for second-order Riesz transforms (in the classical setting). We will focus on the weak-type bounds for $1<p<2$ - the remaining estimates can be treated analogously. Suppose that $\lambda_p$ is the best constant in the estimate
$$ \mathbb{E}(|G_\infty|-\lambda_p)_+\leq \E |F_\infty|^p,$$
valid for all pairs  $(F,G)$ of finite martingales starting from $0$ such that $G$ is a $\pm 1$-transform of $F$. The value of $\lambda_p$ appears in the statement of Theorem 9 above, the fact that it is already the best for martingale transforms follows from the examples exhibited in \cite{O4}. For any $\e>0$, Theorem \ref{laminate} yields the existence of $f:\R^2\to \R$, supported on the unit disc, such that
$$ \int_{\R^2}\big(|(R_1^2-R_2^2)f|-\lambda_p+\e\big)_+\mbox{d}x\mbox{d}y>\int_{\R^2}|f|^p\mbox{d}x\mbox{d}y.$$
That is, if we set $A=\{|(R_1^2-R_2^2)f|\geq \lambda_p-\e\}$, we get
\begin{equation}\label{lower_bound_Riesz}
 \int_A |(R_1^2-R_2^2)f|\mbox{d}x\mbox{d}y>\int_{\R^2}|f|^p\mbox{d}x\mbox{d}y+(\lambda_p-\e)|A|.
\end{equation}
However, if the weak-type estimate holds with a constant $c_p$, Young's inequality implies
\begin{align*}
 \int_A|(R_1^2-R_2^2)f|\mbox{d}x\mbox{d}y&\leq \int_{\R^2}|f|^p\mbox{d}x\mbox{d}y+\frac{(p-1)c_p^{p/(p-1)}}{p^{p/(p-1)}}|A|.
 \end{align*}
Therefore, the inequality \eqref{lower_bound_Riesz} enforces that
 $$ \frac{(p-1)c_p^{p/(p-1)}}{p^{p/(p-1)}}\geq \lambda_p$$
 (since $\e$ was arbitrary). This estimate is equivalent to 
$$c_p\geq\left(\frac{1}{2}\Gamma\left(\frac{2p-1}{p-1}\right)\right)^{1-1/p},$$
which is the desired sharpness.

\subsection{From continuous to discrete sharp estimates}

We claim that the sharp bounds found for the continuous second order Riesz
transforms also hold in the case of purely discrete groups. Groups of mixed
type would be treated in the same manner.
We illustrate those results only for the sharpness in Theorem \ref{T: LpRieszsemidiscrete}
and in Theorem \ref{T: weak-type estimate} since other results follow the same lines.
Precisely, we show that the sharpness in the
discrete case is inherited from the sharpness of the continuous case through
the use of the so--called fundamental theorem of finite difference methods
from Lax and Richtmyer {\cite{LaxRic1956a}} (see also {\cite{LeV2007a}}). This
result states that {\tmstrong{stability}} and {\tmstrong{consistency}} of the
finite difference scheme implies {\tmstrong{convergence}} of the approximate
finite difference solution towards the continuous solution, in a sense that we
detail below.

\paragraph{Finite difference Riesz transforms}Let $u=R^{2}_{i} f$ be the
$i$-th second order Riesz transform in $\Omega \assign \mathbbm{R}^{N}$ of a
function $f \in L^{p}$. The function $u$ is the unique solution to the Poisson
problem in $\mathbbm{R}^{N}$, $\Delta u= \partial^{2}_{i} f$ in
$\mathbbm{R}^{N}$ (see {\cite{Eva1998}}). This is a problem of the form $A u =
B f$, where $A= \Delta$ and $B= \partial^{2}_{i j}$. Introduce now a finite
difference grid of step--size $h>0$, that is the grid $\Omega_{h} \assign h
\mathbbm{Z}^{N}$. The functions $v_{h}$ defined on $\Omega_{h}$ are equipped
with the $L^{p}_{h}$ norm defined as
\[ \| v_{h} \|_{L_{h}^{p}}^{p} \assign \sum_{x \in \Omega_{h}} | v_{h} ( x )
   |^{p}  h^{N}. \]
It is common to identify a finite difference function $v_{h}$ defined on the
grid $\Omega_{h}$ with the piecewise constant function (also denoted) $v_{h}
:\mathbbm{R}^{N} \rightarrow \mathbbm{C}$ such that $v_{h} ( x ) =v_{h} ( y )$
for all $x$'s in the open cube $\Omega ( y )$ of volume $h^{N}$ centered
around the grid point $y \in \Omega_{h}$. With this notation, we might write
finite difference integrals in the form
\[ \| v_{h} \|_{L_{h}^{p}} = \int_{x \in \mathbbm{R}^{N}} | v_{h} ( x ) |
   \mathd \mu_{h} ( x ). \]
The finite difference second order Riesz transform $u_{h} =R^{2}_{i} f_{h}$ of
$f_{h}$ is the solution to the problem $A_{h} u_{h} =B_{h}^{2} f_{h}$, where
$A_{h} \assign \Delta_{h}$ is the finite difference Laplacian and $B_{h}
\assign \partial_{i,h}^{2}$ the $3$--point finite difference second order
derivative. Precisely, for any $x \in \Omega_{h}$, any $v_{h} : \Omega_{h}
\rightarrow \mathbbm{R}$,
\[ ( \partial_{i,h}^{2} v_{h} ) ( x ) \assign \frac{v_{h} ( x+h e_{i} )
   -2v_{h} ( x ) +v_{h} ( x-h e_{i} )}{h^{2}} \]
\[ ( \Delta_{h} v_{h} ) ( x ) \assign \sum_{i=1}^{N} ( \partial_{i,h}^{2}
   v_{h} ) ( x ). \]
It is classical that we have the {\tmstrong{consistency}} of the discrete
problem with respect to the continuous problem, that is for given smooth
functions $u$ and $f$ we have $\Delta_{h} u= \Delta u+ \mathcal{O}( h )$ and
$\partial^{2}_{i,h} f= \partial_{i}^{2} f+ \mathcal{O}( h )$, where the
coefficients in $\mathcal{O}( h )$ include as a factor up to fourth--order
derivatives of $u$ or $f$. This implies in particular that $B_{h} f=B f+
\mathcal{O} ( h )$ in $L_{h}^{p}$ for any given smooth function $f$ with
compact support. It is also classical that $( - \Delta_{h} )^{-1}$ is bounded
in $L_{h}^{p}$ uniformly w.r.t. $h$. This is the $L^{p}$
{\tmstrong{stability}} of the finite difference scheme. The fundamental
theorem of finite difference methods implies the $L^{p}$
{\tmstrong{convergence}} of the sequence of discrete second order Riesz
transforms $u_{h}$ towards the continuous second order Riesz transform $u$.

\paragraph{Discrete Riesz tranforms on Lie-Group}Observe that the finite
difference Riesz transform $u_{h} =R^{2}_{i,h} f_{h}$ defined on the grid
$\Omega_{h}$, also gives rise to a Riesz transform on the Lie group
$\Omega_{1} =\mathbbm{Z}^{N}$. This is a consequence of the homogeneity of
order zero of the Riesz transforms. Indeed, the equation $\Delta_{h} u_{h} =
\partial^{2}_{i,h} f_{h}$ rewrites as $\Delta_{1} u_{1} = \partial^{2}_{i,1}
f_{1}$, where $u_{1} ( y ) \assign u_{h} ( y/h )$, $f_{1} ( y ) \assign f_{h}
( y/h )$ for all $y \in \mathbbm{Z}^{N}$, and where $\Delta_{1}$ and
$\partial^{2}_{i,1}$ are the discrete differential operators defined on
$\mathbbm{Z}^{N}$. We have also $\| u_{h} \|_{L^{p}_{h}} =h^{N/p} \| u_{1}
\|_{L^{p}_{1}}$ and $\| f_{h} \|_{L^{p}_{h}} =h^{N/p} \| f_{1}
\|_{L^{p}_{1}}$. Notice that for all $h$, this ensures that $\| u_{h}
\|_{L^{p}_{h}} / \| f_{h} \|_{L^{p}_{h}} = \| u_{1} \|_{L^{p}_{1}} / \| f_{1}
\|_{L^{p}_{1}}$.

\paragraph{Sharpness for Theorem \ref{T: LpRieszsemidiscrete} in the discrete setting}
In the continuous setting, the sharpness was proved in
\cite{GeiMonSak2010a} based on the combintation $R^2_\alpha=R_1^2-R_2^2$ of second order Riesz tranforms.
Let $u^{( k )} =R_{\alpha}^{2} f^{( k )}$ a sequence of second order Riesz
transforms yielding the sharp constant $C_{p}$ in the estimate, that is $\|
u^{( k )} \|_{p}  /  \| f^{( k )} \|_{p}   \rightarrow  C_{p}$ as $k$ goes to
infinity. For each $k \in \mathbbm{N}$ and $h>0$, introduce the finite
difference approximation $f_{h}^{( k )}$ of $f^{( k )}$ and the corresponding
finite difference Riesz transform. Thanks to the convergence of the finite
difference scheme, we can extract a subsequence $f_{h_{k}}^{( k )}$ such that
$\| u_{1}^{( k )} \|_{p}  /  \| f_{1}^{( k )} \|_{p}  = \| u_{h_{k}}^{( k )}
\|_{p}  /  \| f_{h_{k}}^{( k )} \|_{p}   \rightarrow  C_{p}$. Therefore
$C_{p}$ is also the sharp constant for the second order Riesz transforms in
$\mathbbm{Z}^{N}$.

\paragraph{Sharpness for Theorem \ref{T: weak-type estimate} in the discrete setting}
Recall that we have a bound of the form
\[ \| R^{2}_{\alpha} f \|_{L^{p, \infty} ( \mathbbm{G},\mathbbm{C} )} \assign
   \sup_{E} \left\{ \mu_{z} ( E )^{1/p-1} \int_{E} | R^{2}_{\alpha} f | \mathd
   \mu_{z} \right\} \leqslant C_{p} \| f \|_{L^{p}} \]
for a certain constant $C_{p}$ that is known to be sharp in the case of
continuous second order Riesz transforms. In order to prove sharpness when the
Lie group $\mathbbm{G}$ does not have enough continuous components, it suffices again
to approximate a sequence of continuous extremizers by a sequence of finite
difference approximations. Take $\mathbbm{G}=\mathbbm{R}^{N}$. For any
$\varepsilon >0$, let $f$, $u \assign R^{2}_{\alpha} f$, and $E$ with finite
measure chosen so that
\[ \mu_{z} ( E )^{1/p-1}  \| u \|_{L^1(E)}   / \| f \|_{L^{p}} \geqslant C_{p} - \varepsilon. \]
We can assume without loss of generality that $f$ is a smooth function with compact support.
Let $f_{h}$ a finite difference approximation of $f$ defined as its $L^2$ projection on the grid, and $u_{h}$ its discrete
second order Riesz transform both defined on $\Omega_{h} \assign
h\mathbbm{Z}^{N}$. Since $\mu_{z} ( E )$ is the finite $N$--dimensional Lebesgue
measure of $E$, we use outer measure approximations of $E$ followed by approximations from below
by a finite number of small enough cubes of size $h$ centered around the grid points of $\Omega_h$,
to define a ``finite difference" approximation $E_h$ of $E$ such that
\[
\mu_{h} ( E_h ) \assign \sum_{x\in E_h} h^N \rightarrow \mu_{z} ( E )
\]
when $h$ goes to zero. Since the discrete Riesz transforms are stable in $L^2$,
the Lax-Richtmyer theorem ensures that $\| u_h \|_{L^2_h}  \rightarrow \| u \|_{L^2}$
which implies $\| u_h \|_{L^1_h(E)}  \rightarrow \| u \|_{L^1(E)}$
and also $\| u_h \|_{L^1_h(E_h)}  \rightarrow \| u \|_{L^1(E)}$.
Therefore for $h$ small enough,

\[ _{} \mu_{h} ( E_h )^{1/p-1} \| u_h \|_{L^1_h(E_h)}
   / \| f_{h} \|_{L^{p}} \geqslant C_{p} -2 \varepsilon. \]
Let as before $u_{1} ( y ) \assign u_{h} ( y/h )$, $f_{1} ( y ) \assign f_{h}
( y/h )$ for all $y \in \Omega_{1} \assign \mathbbm{Z}^{N}$, and $E_{1}
\assign E/h$. We have sucessively $\mu_{h} ( E ) =h^{N} \mu_{1} ( E_{1} )$,
$ \| u_h \|_{L^1_h(E_h)}  =  h^{N} \| u_1 \|_{L^1_h(E_1)} $
and $\| f_{h} \|_{L^{p}} =h^{N/p} \| f_{1} \|_{L^{p}}$. This yields immediately
\[ _{} \mu_{1} ( E_{1} )^{1/p-1} \| u_1 \|_{L^1_h(E_1)}   / \| f_{1} \|_{L^{p}} \geqslant C_{p} - 2 \varepsilon , \]
allowing us to prove sharpness for the class of discrete groups we are
interested in.

\end{document}